\documentclass{amsart}

\usepackage{amsmath, amssymb, amsfonts, amsthm}
\usepackage{xspace}
\usepackage{mathtools}
\usepackage[all, 2cell]{xy}
\usepackage{tikz}
\usepackage{tikz-cd}
\usetikzlibrary{matrix}
\usetikzlibrary{arrows}
\usetikzlibrary{positioning}
\usetikzlibrary{decorations.pathmorphing}

\theoremstyle{plain}
\newtheorem{main}{Theorem}
\newtheorem{mcorollary}[main]{Corollary}
\newtheorem{theorem}{Theorem}[section]
\newtheorem{lemma}[theorem]{Lemma}
\newtheorem{corollary}[theorem]{Corollary}
\newtheorem{proposition}[theorem]{Proposition}

\newcommand{\bR}{\mathbb{R}}
\newcommand{\bQ}{\mathbb{Q}}
\newcommand{\bC}{\mathbb{C}}
\newcommand{\bZ}{\mathbb{Z}}
\newcommand{\bA}{\mathbb{A}}

\newcommand{\fF}{\mathfrak{f}}
\newcommand{\fP}{\mathfrak{p}}
\newcommand{\fQ}{\mathfrak{Q}}

\newcommand{\cO}{\mathcal{O}}

\newcommand{\cA}{\mathcal{A}}

\newcommand{\MED}{\medskip \hspace{2ex} \medskip}

\newcommand{\GL}{{\rm GL}}

\begin{document}

\title{Multiplicity One Theorems, $S$-version}
\date{}
\author{Song Wang}
\address{Song Wang, Academy of Mathematics
 and Systematics Sciences, the Morningside Center of Mathematics,
 the Key Laboratory of Hua, Chinese Academy of Science.}
\thanks{The author was supported in part by national 973 project 2013CB834202,
also National Natural Science Foundation of China (NNSFC 11321101), and also by the
One Hundred Talent's Program from Chinese Academy of Science.}
\email{\texttt{songw1973@amss.ac.cn}}
\maketitle

\section{Introduction} \label{S:1}


Let $K$ be a number field, and $\pi$ a cuspidal automorphic representation
of \\
$\GL_{d} (\bA_{K})$. It is well known that $\pi$ is unique determined by the
Satake parameter $c (\pi, v)$ for almost all $v$, and even more precisely,
almost all $v$ of degree $1$. Also, it suffices for us to test only finitely
many $v$, namely, there is a bound $\mathfrak{N}$ which is some expression
in terms of $K$, $d$, $N (\pi)$ such that there is a $v$ with
$\pi_{v} \cong \pi'_{v}$ and $N \fP_{v} < \mathfrak{N}$ (\cite{Brumley06},
\cite{Li2010}, \cite{L-M-O79}, \cite{L-O77}, \cite{LW2009}, \cite{Moreno85},
\cite{Wyh2006}).

\medskip

In this article, we are going to prove some $S$-effective version
of multiplicity one theorems. Roughly speaking, all situation as above,
if $\pi$ and $\pi'$ are not equivalent, then there us also a bound $\mathfrak{N} (S)$
which is some expression in terms of $K$, $d$, ${\rm max} (N (\pi), N (\pi')$,
such that there is a $v \notin S$ with
$\pi_{v} \cong \pi'_{v}$ and $N \fP_{v} < \mathfrak{N}$.
We encounter such $S$-versions in a series work (\cite{Wang2001}, \cite{Wang2013-1}).

\medskip

Now we state our first two theorems. One is the $S$-effective version
of the Chebotarev Density Theorem, and the other one is the multiplicity one
for $GL (1)$. Thrughout,
for any place $v$ of $K$, denote $\fP_{v}$ be the formal prime corresponding to $v$
which denotes also the prime idea of $\cO_{K}$ or $\cO_{K_{v}}$.
$d_{K}$, $d_{L}$ denote the discriminants of $K$ and $L$ respectively,
and $\zeta_{K} (s), \zeta_{L} (s)$ denote the Dedekind zeta functions of $K$ and $L$
respectively. if $\chi$ is a global character of $\bA_{K}^{\times}$,
$N (\chi)$ denotes the conductor of $\chi$ (for explicit definitions see next section.
Also \cite{H-R95}, \cite{Moreno85}, \cite{Ra99}, \cite{Ra2000}, \cite{RW2003}).
For a finite set $S$ of places of $K$, $N_{S}$ denotes the product of norms of
the places in $S$.

\medskip

\begin{main} \label{TM:A}
\textnormal{\textbf{($S$-effective Version of the Chebotarev Density Theorem)}}

\medskip

Let $L / K$ be a finite Galois extension of number fields and $L \ne \bQ$.
Then there is an effectively computable absolute constant $C$ satisfying the following:
For each finite set $S$ of places of $K$, a conjugacy class $[\sigma]$
in ${\rm Gal} (L / K)$, there is a place $v$ of $K$ such that

\medskip

\textnormal{(1)} $v \notin S$ and is of degree $1$.

\medskip

\textnormal{(2)} The Artin symbol $\left( \frac{L / K}{\fP_{v}} \right) = [\sigma]$.

\medskip

\textnormal{(3)}
\[
N (\fP_{v}) \leq {(d_{L} N_{S}^{[L : K]})}^{C}
\]
\end{main}

\medskip

\begin{main} \label{TM:B}
\textnormal{\textbf{($S$-effective Version of the Multiplicity One for $GL (1)$)}}

\medskip

Let $K$ be a number fields and $\chi$ a global character of $\bA_{K}^{\times}$.
Then there is an effectively computable absolute constant $C$ satisfying the following:
For each finite set $S$ of places of $K$, there is a place $v$ of $K$ such that

\medskip

\textnormal{(1)} $v \notin S$ and is of degree $1$.

\medskip

\textnormal{(2)} $\chi_{v}$ is unramified and nontrivial.

\medskip

\textnormal{(3)}
\[
N (\fP_{v}) \leq {(d_{K} N (\chi) N_{S})}^{C}
\]
\end{main}

\medskip

Now we make some definitions.
Let $\pi$ be a unitary cuspidal automorphic representation of
$\GL_{d} (\bA_{K})$, and
\[
L (s, \pi_{\infty}) =
\pi^{-s d n_{K} / 2} \prod_{j = 1}^{d n_{K}} \Gamma ((s + b_{j} (\pi)) / 2)
\]
Moreover, $N (\pi)$ be the level of $\pi$ so that
\[
\Lambda (s, \pi) = L (s, \pi_{\infty}) L (s, \pi)
\]
satisfying
\[
\Lambda (s, \pi) = W (\pi) N (\pi)^{1/2 - s} \Lambda (1-s, \tilde{\pi})
\]
where $W (\pi)$ the root number of $\pi$.

\medskip

Now define the \emph{extended analytic conductor} $C (\pi)$ (in the sense of
\cite{H-R95}, \cite{RW2003} etc.\ )
be $N (\pi) \prod_{j = 1}^{d n_{K}} (1 + |b_{j} (\pi)|)$.

\medskip

Let $\pi'$ be another unitary cuspidal automorphic representation of
$\GL_{d'} (\bA_{K})$, and
\[
L (s, \pi_{\infty} \times \pi'_{\infty}) = \pi^{-s d d' n_{K} / 2}
 \prod_{j = 1}^{d d' n_{K}} \Gamma ((s + b_{j} (\pi \times \pi')) / 2)
\]
Moreover, $N (\pi \times \pi')$ be the level of the Rankin-Selberg
$\pi \times \pi'$ so that
\[
\Lambda (s, \pi \times \pi') = L (s, \pi_{\infty} \times \pi'_{\infty})
 L (s, \pi_{\rm f} \times \pi'_{\rm f})
\]
satisfying
\[
\Lambda (s, \pi \times \pi') = W (\pi \times \pi') N (\pi \times \pi')^{1/2 - s}
 \Lambda (1-s, \tilde{\pi} \times \tilde{\pi}')
\]
where $W (\pi \times \pi')$ the root number of $\pi \times \pi'$.

\medskip

Now define the \emph{extended analytic conductor} $C (\pi \times \pi')$ (in the sense of
\cite{H-R95}, \cite{RW2003}, \ )
be $N (\pi) \prod_{j = 1}^{d d' n_{K}} (1 + |b_{j} (\pi \times \pi')|)$.

\medskip

Moreover, Define the \emph{bound for Ramanujuan} for $\pi$ be the upper log bounds of coefficients
of the cusp form. Namely, for each $\pi \in \cA_{0} (K)$,
its bound for Ramanujuan
${\rm RB} (\pi) = {\rm Sup}_{v} {\rm Max_{j = 1 \ldots d}} \log_{q_{v}} \alpha_{v, i} (\pi)$
where $v$ runs through all finite places of $K$ where $\pi_{v}$ is unramified
and $\alpha_{v, 1}, \ldots \alpha_{v, d}$
are Satake parameters of $\pi_{v}$.

\medskip

\begin{main} \label{TM:C}
Let $\pi$ and $\pi'$ be two unitary cuspidal automorphic representations
of $\GL_{d} (K)$. Let $S$ be a finite set of places of $K$,
and $Q = {\rm max} (C (\pi), C (\pi'))$ and assume that the bound for Ramanujuan for
$\pi$ and $\pi'$ are $< R$.

\medskip

Then if $\pi \not\cong \pi'$,
there exists a place $v$ of $K$ such that $\pi_{v} \not\cong \pi'_{v}$
and
\[
N (\fP_{v}) \leq \begin{cases}
C Q^{1 + \epsilon} N_{S}^{\epsilon} &(d = 1) \notag \\
C Q^{2 d + \frac{d (d - 2)}{d H + 1} + \epsilon}
    N_{S}^{\frac{d^{3} (2 R + H)}{d H + 1} + \epsilon}
&(\text{general $d$}) \notag
\end{cases}
\]
where $C$ is some effectively computable constant
only depending on arbitrarily chosen number $H > 2 R, \epsilon > 0$, $K$ and $d$.
\end{main}

\medskip

\emph{Remark}: When $d \geq 2$, the two ends of $H > 2 R$ lead to the following two estimates:
\[
N (\fP_{v}) <<_{\epsilon} {\rm min} (Q^{2 d + \frac{d (d - 2)}{1 + 2 R H}}
 N_{S}^{d^{2} \frac{4 R d^{3}}{2 d R + 1} \epsilon},
 Q^{2 d + \epsilon} N_{S}^{d^{2} + \epsilon})
\]
In particular, if $R = 0$ (i.e., Ramanujuan holds), then
the two ends are
\[
N (\fP_{v}) <<_{\epsilon} {\rm min} (Q^{2 d + \epsilon} N_{S}^{d^{2} + \epsilon},
Q^{d^{2} + \epsilon} N_{S}^{\epsilon})
\]

\medskip

This theorem is an $S$-effective refinement of \cite{LW2009}. When take $H = 0$,
and be aware that the unitary Grossencharacters are cusp forms of $GL (1)$
with the bound for Ramanujuan as $0$,
we have the following corollary (Theorem 4.1 (B), \cite{Wang2013-1}).

\medskip

\begin{mcorollary}
Let $K$ be a number field and $\chi$ a nontrivial unitary character of $C_{K}$.
Then there is a place $v$ of $K$ such that

\medskip

\textnormal{(1)} $\fP_{v} \notin S$.

\medskip

\textnormal{(2)} $\chi_{v} \ne 1$ and is not ramified.

\medskip

\textnormal{(3 B)} $N (\fP_{v}) <<_{\epsilon, K} N (\chi)^{1/2 + \epsilon} N_{S}^{\epsilon}$
for every $\epsilon > 0$.

\medskip

where $A (\chi, S) = d_{K} N (\chi) N_{S}$.

\end{mcorollary}

\medskip

Now we summarize the technique we used to proved the theorem. Theorem ~\ref{TM:A}
and Theorem ~\ref{TM:B} use a refinement of arguments of
\cite{L-M-O79} and Theorem ~\ref{TM:C} uses Landau's idea plus arguments with modification
of series papers (\cite{Brumley06}, \cite{LW2009}, \cite{Wyh2006}).
Results in \cite{Wang2013-1} are just the special case.

\medskip

This paper is a byproduct of my projects on the effective version of the Grunwald-Wang
and here I express thanks to my advisor Dinakar Ramakrishnan for the introduction
of the problem and the guidance and help during my student
year and continuing years afterwards.
Finally, the author acknowledges the support in part by national 973 project 2013CB834202,
also National Natural Science Foundation of China (NNSFC 11321101), and also by the
One Hundred Talent's Program from Chinese Academy of Science.

\bigskip

\section{Notations and Preliminaries} \label{S:2}

In this section, we recall certain notations and the preliminaries to be used in the proofs.
Lots of standard results can be found in various textbooks. Experts can skip most parts of
this section.

\medskip

\subsection*{Hecke Characters and Hecke $L$-functions} \MED


First, we recall some basic concepts and facts. Throughout,
$K$ denotes a number field and $\bA_{K}$, $\bA^{\times}_{K}$
denote the Adele rings and the idele group of $K$. It is well
known that $K$, $K^{\times}$ embed into $\bA_{K}$ and $\bA_{K}^{\times}$
respectively.

\medskip

A Hecke character $\chi$ of $K$ has two types of expression:
it is a character of the group fractional ideals $J_{K}$ of $K$,
or a continuous character of $\bA_{K}^{\times} / K^{\times}$. $\chi$
is also called a \emph{global character}.

\medskip

Let $K_{v}$ denote the completion of $K$ at a place $v$ with
$\cO_{K_{v}}$, $\fP_{v}$, $q_{v}$ the valuation ring, the prime ideal
and the residue size $\# (\cO_{K_{v}} / \fP_{v})$ of $K_{v}$ respectively. A continuous
character of $K_{v}$ is called a \emph{local character}.
Given a global character $\chi$, let $\chi_{v}$ be the restriction
of $\chi$ to $K_{v}^{\times}$, then $\chi (x) = \prod_{v} \chi_{v} (x_{v})$
if $x = (x_{v}) \in \bA_{K}^{\times}$ and the product is finite for each $x$
since all but finitely many $\chi_{v}$ are unramified.

\medskip

Now define the conductor of $\chi_{v}$ and $\chi$. For each finite place
$v$, we define the \emph{arithmetic conductor} $N (\chi_{v})$ of $\chi_{v}$
as
\[
N (\chi_{v}) = \begin{cases}
1 &\text{(if $\chi_{v}$ is unramified, i.e., $\chi_{v} (\cO_{K_{v}}) = 1$)} \\
q_{v}^{r} &\text{(if $r$ is the smallest
    integer such that $\chi_{v} (1 + \fP_{v}^{r}) = 1$)}
\end{cases}
\]
and $N (\chi) = \prod_{v} N (\chi_{v})$.

\medskip

Recall that the $L$-function $L (s, \chi) = \prod_{v < \infty} L (s, \chi_{v})$
and the complete $L$-function $\Lambda (s, \chi) = L_{\infty} (s, \chi) L (s, \chi)
= \prod_{v} L (s, \chi_{v})$ where the Gamma factor
is $L_{\infty} (x, \chi) = \prod_{v | \infty} L (s, \chi_{v})$. Here the local $L$-factors
are
\begin{align}
L &(s, \chi_{v}) &&\notag\\
&=  \pi^{- (s + d) / 2} \Gamma ((s + d) / 2)
&&(\text{$v$ is real, $\chi_{v} (x) = {|x|}^{d}$ or ${|x|}^{d - 1} {\rm Sgn}$}) \notag \\
&= 2 {(2 \pi)}^{- (s + d)} \Gamma (s + d)
&&(\text{$v$ is complex, $\chi_{v} (z) = z^{a} \bar{z}^{b}$ with $d = {\rm max} (a, b)$})
\notag \\
&= \frac{1}{1 - \chi_{v} (\mathcal{\pi}_{v}) q_{v}^{-s}}
 &&(\text{$v$ is finite and $\chi_{v}$ is unramified}) \notag \\
&= 1  &&(\text{$v$ is finite and $\chi_{v}$ is ramified}) \notag
\end{align}
where $\mathcal{\pi}_{v}$ is a uniformizer of $K_{v}$, i.e.,
a generator of the unique prime ideal
(also denoted $\fP_{v}$) of $\cO_{K_{v}}$.

\medskip

When $\chi = 1$ then $L (s, \chi) = \zeta_{K} (s)$. Moreover, it is well known
that $L (s, \chi)$ is entire of order $1$ if $\chi \ne 1$ is finite order. Moreover,
$L (s, \chi)$ satisfies the following functional equation:
\[
\Lambda (s, \chi) = W (\chi) {A (\chi)}^{1/2 - s} \Lambda (1 - s, \chi^{-1})
\]
if $\chi$ is unitary. Here $A (\chi) = d_{K} N (\chi)$ is the \emph{analytic conductor}
of $\chi$.

\medskip

Let us list out some facts to be used in later parts of this paper.
One is the Euler expression:
\[
- \frac{L' (s, \chi)}{L (s, \chi)}
= \sum_{v < \infty, \chi \text{is unramified}}
 \sum_{k = 1}^{\infty} \chi_{v}^{k} (\pi_{v}) q_{v}^{-n s} \log q_{v}
\]
which can be easily deduced from the Euler product.
Moreover, by the class field theory, if $\chi$ is of finite order,
then $\chi$ is associated to a cyclic extension
$L / K$, and in fact, it is associated to a $1$-dimensional representation
of ${\rm Gal} (L / K)$ (see the preliminary on Galois representation).
Moreover, Weil (\cite{Weil52}, \cite{Weil55}, \cite{CF79})
showed that any Hecke character $\chi$ is an unramified twist of a character
of finite order unless it is a CM character.

\medskip

\subsection*{Cuspidal and Isobaric Automorphic
 Representations and Automorphic $L$-functions} \MED

Now we recall some notations and facts the automorphic forms
on $\GL (d)$. Most details can be found in various literatures
(\cite{Bump}, \cite{Ra99} etc).

\medskip

Denote $\cA (d, K)$ (resp.\ $\cA_{0} (d, K)$) the set of unitary irreducible
automorphic (resp.\ cuspidal automorphic) representations of $\GL_{d} (\bA_{K})$
and $\cA (K)$ (resp.\ $\cA_{0} (K)$) the set of unitary irreducible
automorphic (resp.\ cuspidal automorphic) representations of $\GL_{d} (\bA_{K})$
for some $d$. Denote $L (s, \pi) = \prod_{v < \infty} L (s, \pi_{v})$ be the (finite part)
$L$-function associated to $\pi$, $L_{\infty} (s, \pi) = \prod_{v | \infty} L (s, \pi_{v})$,
$\Lambda (s, \pi) = L (s, \pi) L_{\infty} (s, \pi)$ be
the infinite part, and the complete $L$-function
of $\pi$ respectively.

\medskip

For $\pi \in \cA (K)$,
It is well known that $\pi = \otimes \pi_{v}$ where $\pi_{v}$ is an irreducible
admissible representation of $\GL_{d} (K_{v})$ for each $v$.
Moreover, by the Langlands classification, for each $\pi \in \cA (K)$, there
are $\pi_{1}, \ldots, \pi_{r} \in \cA_{0} (K)$ such that we have the
following isobaric sum decomposition (\cite{JS}, \cite{JPSS})
$\pi = n_{1} \pi_{1} \boxplus \ldots \boxplus n_{r} \pi_{m}$.
In particular, $L (s, \pi) = \prod_{i = 1}^{r} {L (s, \pi_{i})}^{n_{i}}$,
$L_{\infty} (s, \pi) = \prod_{i = 1}^{r} {L_{\infty} (s, \pi_{i})}^{n_{i}}$.

\medskip

For $\pi \in \cA_{0} (K)$, it is also well known that, $L (s, \pi)$
and $\Lambda (s, \pi)$ are
entire unless $\pi$ is a trivial character so that $L (s, \pi) = \zeta_{K} (S)$ has
a simple pole at $s = 1$.
In fact, we have $\Lambda (s, \pi) = W (\pi) {Q (\pi)}^{1/2 - s}
\Lambda (1 - s, \tilde{\pi})$
where $\tilde{\pi}$ is the contragredient of $\pi$. Here $N (\pi)$ is some positive
integer called the \emph{level} or the \emph{analytic conductor}
of $\pi$, and it is closedly related to the \emph{arithmetic conductor}
which is defined in an arithmetic way. In particular, if $\pi \in \cA_{0} (1, K)$,
then $\pi$ is in fact a unitary continuous idele class character $\chi$
of $\bA_{K}^{\times} / K^{\times}$, and the arithmetic conductor
of $\chi$ is $N (\chi)$ and the analytic conductor is $A (\chi) = d_{K} N (\chi)$.
Moreover $W (\pi)$ is a constant of complex number of absolute value $1$.

\medskip

Now we talk about the local $L$-factors. Let $\pi \in \cA_{0} (d, K)$.
For almost all finite place $v$, $\pi_{v}$ is unramified. In fact, $\pi_{v}$
is a constituent of the unramified (normalized) induced representation
${\rm Ind}_{P}^{\GL_{n}} (\chi_{1} \otimes \ldots \otimes \chi_{d})$
where $\chi_{1}, \ldots, \ldots \chi_{d}$ are unramified characters of
$K_{v}^{\times}$. Let $a_{v, i} = \chi_{i} (\pi_{p})$. Then we have
\[
L (s, \pi_{v}) = \prod_{i = 1}^{d} {(1 - a_{v, i} q_{v}^{-s})}^{-1} =
= {{\rm Det} (1 = q^{-s} c (\pi, v))}^{-1}
\]
where \emph{Satake Parameter} $c (\pi, v)$ is the conjugacy class
of ${\rm diag} (a_{v, 1}, \ldots, a_{v, d})$.
Moreover, if $v$ is an infinite place, then
\[
L (s, \pi_{v}) = \prod_{i = 1}^{d} \Gamma_{v} (s + b_{v, j} (\pi))
\]
for some $b_{v, j} (\pi) (j = 1, \ldots d)$.
Here $\Gamma_{v} (s)$ is defined as $\Gamma_{\bR} (s) = \pi^{-s/2} (s / 2)$
and $\Gamma_{\bC} (s) = 2 (2 \pi)^{-s} \Gamma (s)$.
$b_{v, j} (\pi)$ are called the \emph{infinite type constants}.
Put $L (s, \pi_{v})$ for $v | \infty$ together, we get $L_{\infty} (s, \pi)$
and also the infinite type constant $b_{j} (\pi)$ $j = 1 \dots d n_{K}$.
In fact, $\{b_{j} (\pi)\}$ consists
of those $b_{v, j} (\pi)$ for $v$ real and $b_{v, j} (\pi)$ and $b_{v, j} (\pi) + 1$
for $v$ complex.

\medskip

Now let's talk about the ``extended analytic conductor'' $C (\pi)$.
It is widely used in estimation. Unlike $N (\pi)$, it shows up quite differently
in different literature. However, it is essentially the same for the purpose
of the estimations. In fact, as
$L_{\infty} (s, \pi) = \pi^{- d n_{K} / 2} \prod_{j = 1}^{d n_{K}} \Gamma ((s + b_{j} (\pi)) / 2)$
we define
\[
C (\pi) = N (\pi) \prod_{i = 1}^{d n_{K}} (1 + |b_{j} (\pi)|)
\]
In general, (at least for all cases we know) ${\rm Re} b_{j} (\pi) \geq 0$.

\medskip

Now look at $\tilde{\pi}$. Note that $N (\tilde{\pi}) = N (\pi)$,
As $\pi$ is unitary, $\tilde{\pi} = \bar{\pi}$ so that
all related quantities associated to $\tilde{\pi}$ are complex conjugate of those of
$\pi$. In particular, the Satake parameter $c (\tilde{\pi}, v) = \overline{c (\pi, v)}$
when $\pi_{v}$ is unramified.

\medskip

Later on, we will see that $L (s, \pi)$ is standard.

\medskip

\subsection*{Galois representations and Artin $L$-functions} \MED

By the class field theory, when $\chi$ is of finite order $m$,
then
$\chi$ is associated to a cyclic Galois extension $L / K$ of degree $m$
and $L (s, \chi) = L (s, \rho)$ where $\rho$ is the $1$ dimensional
complex representation of ${\rm Gal} (L / K)$ associated to $\chi$.
In particular, $v$ splits in $L / K$ if and only if $\chi_{v}$
is trivial. So this is an example of Artin $L$-functions.

\medskip

In general, let $\rho$ be a Galois representation, namely,
a finite dimensional representation of ${\rm Gal} (L / K)$
where $L$ is a finite Galois extension of $K$. Recall that we can
define the $L$-function as the following: Let $V$ be the representation space
of $\rho$ and $n = {\rm dim}_{\bC} V$.
For each finite place $v$ of $K$, $w$ a place of $L$ over $v$
let $\sigma_{w / v} \in {\rm Gal} (L /K)$ be of the Frobenious element of $L / K$
at $w / v$, i.e., $\sigma_{w / v} \in {\rm Gal} (L / K)$ such that (i)
$\sigma_{w / v} (\fP_{w}) = \fP_{w}$
where $\fP_{w}$ the prime idea of $L$ associated to $w$,
(ii) $\sigma_{w / v} (x) \equiv x^{q_{w}}$
 where $q_{w}$ be the size of $\cO_{L_{w}} / \fP_{w}$.
It is well known that, the conjugacy class of $\sigma_{w / v}$ is independent
of the choice of $w$. Now define
\[
L (s, \rho_{v}) = {{\rm det}_{V^{I_{w / v}}} (I - q_{v}^{-s} \sigma_{w / v})}^{-1}
\]
where $I_{w / v}$ be the inertial group of $L / K$ at $w / v$, and $V^{I_{w / v}}$
be the subspace of $V$ fixed by $I_{w / v}$. It is easy to see
that $L (s, \rho_{v})$ is independent of the choice of $w$.

\medskip

Define
\[
L (s, \rho) = \prod_{v} L (s, \rho_{v})
\]
We can also define the Artin conductor
 $A_{L / K} (\rho)$ (the explicit definition is a little bit
complicated, and can be found in a lot of literatures.)
It is well known that (1) $L (s, \rho)$ is meromorphic, (2) it allows a
functional equation:
\[
\Lambda (s, \rho) = \pi^{-n s / 2} {\Gamma (s / 2)}^{a (\rho)}
 {\Gamma ((s + 1) / 2)}^{b (\rho)}
\]

\medskip

Recall that The Artin conjecture asserts that $L (s, \rho)$
is entire if $\rho$ is a nontrivial
irreducible representation. Moreover, the strong Artin conjecture asserts that,
in this case there is a unitary cuspidal automorphic representation
$\pi$ of $\GL_{d} (\bA_{K})$ such that $L (s, \rho) = L (s, \pi)$.
In this case, not only $L (s, \rho)$ has good analytic property
but so do also a large families of $L$-functions, eg.\ the Rankin-Selberg $L$-functions
(See later subsections).
Moreover, by the standard argument which could be seen in a lot of literatures,
the analytic Artin conductor $A_{L / K} (\rho)$ is the same as $N (\pi)$,
the analytic conductor. In particular, this gives rises to internal
relations of arithmetic features (for example, conductors) between
$\rho$ and $\pi$.

\medskip

\subsection*{Rankin-Selberg $L$-functions} \MED

Now we recall the fact of Rankin-Selberg $L$-functions, which provide a large family of
$L$-functions with good analytic properties. (For details, see.)

\medskip

Let $\pi$ and $\pi'$ be two isobaric automorphic representations of
$\GL_{d} (\bA_{K})$ and $\GL_{d'} (\bA_{K})$. Then we can define
the Rankin-Selberg product $L (s, \pi \times \pi')$ which is an
$L$-function of degree $d d' n_{K}$. We know there are various ways to define the
Rankin-Selberg $L$-functions. One way is a formal one, namely, by the local Langlands.
For each place $v$ of $K$,  let $W'_{K_{v}}$ be the
Weil-Delign group of the local field $K_{v}$,
let $\sigma (\pi_{v})$ be the $d$-dimensional admissible representations
of $W'_{K_{v}}$ associated to $\pi_{v}$, and $\sigma (\pi'_{v})$ defined similarly.
Define the local $L$-factors
\[
L (s, \pi_{v} \times \pi'_{v}) = L (s, \sigma (\pi_{v}) \otimes \sigma (\pi'_{v}))
\]
where the local $L$-factors of Weil-Delign group representations are defined similarly
as the $L$-functions of local Galois representations.
In particular, if $v$ is a finite place and $\pi$ and $\pi'$ are unramified at $v$, and
\[
L (s, \pi_{v}) = \prod_{i = 1}^{d} {(1 - a_{v, j} (\pi) q_{v}^{-s})}^{-1}
\]
\[
L (s, \pi'_{v}) = \prod_{i = 1}^{d'} {(1 - a_{v, j} (\pi') q_{v}^{-s})}^{-1}
\]
Then
\[
L (s, \pi_{v} \times \pi'_{v}) = \prod_{i = 1}^{d} \prod_{j = 1}^{d'}
{(1 - a_{v, j} (\pi) a_{v, j} (\pi') q_{v}^{-s})}^{-1}
\]

\medskip

Like automorphic $L$-functions, we can define
$L (s, \pi \times \pi') = \prod_{v < \infty} L (s, \pi_{v} \times \pi'_{v})$,
$L_{\infty} (s, \pi \times \pi') = \prod_{v} L (s, \pi_{v} \times \pi'_{v})$
and $\Lambda (s, \pi \times \pi') = L (s, \pi \times \pi') L_{\infty} (s, \pi \times \pi')$.

\medskip

We can also define these through certain zeta integrals.

\medskip

It is well known that $L (s, \pi \times \pi')$ converges on some right half plane of
$s$, and extends to a meromorphic function.
There is a positive integer
$N = N (\pi \times \pi')$ and a number $W = W (\pi \times \pi')$ of norm $1$
such that
\[
\Lambda (s, \pi \times \pi') = W N^{1/2 - s} \Lambda (1 - s, \tilde{\pi} \times \tilde{\pi'})
\]
Using the language above,
$L (s, \pi \times \pi')$ is standard when $\pi$ and $\pi'$ are isobaric, with its dual
$L (s, \tilde{\pi} \times \tilde{\pi}')$.

\medskip

\begin{proposition} \label{T:201}
Let $\pi$ and $\pi'$ be two isobaric automorphic representations
of \linebreak
$\GL_{d} (\bA_{K})$ and $\GL_{d'} (\bA_{K})$ respectively.

\medskip

\textnormal{(1)} If $\pi$ and $\pi'$ are cuspidal,
then $L (s, \pi \times \pi')$ is entire unless $\pi' \cong \tilde{\pi} \otimes {||}^{i t_{0}}$,
and in this case, $L (s, \pi \otimes \pi')$ has a simple pole at $s = 1 - i t_{0}$, and is holomorphic
else where.

\medskip

\textnormal{(2)} If $\pi = \boxplus m_{i} \pi_{i}$ and $\pi' = \boxplus n_{j} \pi'_{j}$
are the isobaric decompositions of $\pi$ and $\pi'$. Then $L (s, \pi \times \pi')$
is meromorphic, and all of its poles lie on the line ${\rm Re} (s) = 1$. Moreover,
$s = 1 - i t_{0}$ is a pole if and only if $\pi'_{j} \cong \tilde{\pi_{i}} \otimes {||}^{i t_{0}}$
for some $i$ and $j$, and the order of the pole at $s = 1 - i t_{0}$ is $\sum m_{j} n_{k}$
where the sum is taken over $(j, k)$ such that $\pi'_{k} \cong \tilde{\pi}_{j} \otimes {||}^{i t_{0}}$.

\medskip

\textnormal{(3)} In particular, $\pi$ is cuspidal if and only if $L (s, \pi \times \tilde{\pi})$
has a simple pole at $s = 1$.
\end{proposition}

\medskip

\qedsymbol

\medskip

Note that, the automorphic $L$-functions are the special case
of the Rankin-Selberg $L$-functions.

\medskip

Now we define the extended analytic conductor.

Write
$L_{\infty} (s, \pi \times \pi') = \pi^{- d^{2} n_{K} / 2}
    \prod_{j = 1}^{d^{2} n_{K}} \Gamma ((s + b_{j} (\pi \times \pi')) / 2)$
we define
\[
C (\pi \times \pi') = N (\pi \times \pi') \prod_{i = 1}^{d^{2} n_{K}} (1 + |b_{j} (\pi \times \pi')|)
\]
In general, (at least for all cases we know) ${\rm Re} b_{j} (\pi) \geq 0$.

\medskip

We quote the following estimation for (extended) analytic conductor Rankin-Selberg product:

\medskip

\begin{proposition} \label{T:202}
Let $\pi \in \cA (d, K)$ and $\pi' \in \cA (d', K)$.
Then there are constant $C' > 0$ depending on $n_{K}$
such that
\[
N (\pi \times \pi') \leq {N (\pi)}^{d'} {N (\pi')}^{d}
\]
and
\[
C (\pi \times \pi') \leq C' {C (\pi)}^{d'} {C (\pi')}^{d}
\]
\end{proposition}

\medskip

See \cite{He2000}. (Can choose $C' = 2^{d d' n_{K}}$.) \qedsymbol

\medskip

\subsection*{Bounds for Ramanujuan} \MED

In 1916, Ramanujuan conjectured that the Ramanujuan $\Delta$ function
has its coefficients $\tau (n) \leq << n^{11/2 + \epsilon}$. more than half century later,
Deligne proved this conjecture. In fact, the most general formulation of the Ramanujuan
conjecture is formulated as, for each $\pi \in \cA_{0} (d, K)$, and each place $v$
of $K$, the Satake parameters $ a_{v, j}$ are all purely imaginary numbers. (In this case,
$\pi_{v}$ is called \emph{tempered}).

\medskip

Unfortunately, this is still open, even for the case $\GL_{2} (\bQ)$. In fact,
holomorphic modular forms satisfies the Ramanujuan conjecture but for the Maass wave forms,
this is the Selberg conjecture, and is still open.

\medskip

In fact, for each $\pi \in \cA_{0} (K)$, we define the following bound called
\emph{bound for Ramanujuan} for $\pi$
\[
{\rm RB} (\pi) = {\rm Max}_{v, j = 1, \ldots, d} {\rm Re} \log_{q_{v}} a_{v, j} (\pi)
\]
and say $\pi$ satisfies $H_{d} (\delta)$ if for (sufficiently large) $v$,
${\rm Re} \log_{q_{v}} a_{v, j} (\pi) \leq \delta$. See \cite{BB}.
We say that $H_{d} (\delta)$ holds if each $\pi \in \cA_{0} (K)$ satisfies
$H_{d} (\delta)$.

\medskip

\emph{Remark}: We have, by \cite{BB}, $H_{n} (R_{n})$ holds
for $R_{2} = 7/64$, $R_{3} = 5/14$, $R_{4} = 9/22$
and $R_{n} = 1/2 - 2/(n (n + 1) + 2)$ for general $n$.

\medskip

\emph{Remark}: One is easy to see that, if $\pi \in \cA_{0} (K)$ has
$R$ as bound for Ramanujuan. Then $L (1 + H + 2 R, \pi \times \tilde{\pi})
 \leq {(\zeta_{K} (1 + H))}^{d^{2}}$
for $H > 0$.

\medskip

\subsection*{A lemma on coefficients of Rankin-Selberg $L$-functions of positive type}
\MED

We quote a crucial lemma on coefficients of Rankin-Selberg $L$-functions of positive type.
(cf.\ Lemma 2 of \cite{Brumley06}, Proposition 4 of \cite{Wyh2006} and Section 3 of \cite{LW2009}.)

\medskip

\begin{proposition} \label{T:20X}
Let $\pi \in \cA_{0} (d, K)$, and $L (s, \pi \times \tilde{\pi}) = \sum_{n = 1} a_{n} n^{-s}$
be the Dirichlet expansion. For each place $v$ of $K$, we have
\[
a_{q_{v}^{d}} \geq 1
\]
\end{proposition}

\qedsymbol

\medskip

\emph{Remark}: This is a corollary of certain properties of Schur polynomials.

\bigskip

\section{L-O Method} \label{S:3}

In this section, we prove Theorem ~\ref{TM:A} and Theorem ~\ref{TM:B}. We mainly follow the method
of \cite{L-O77} and \cite{L-M-O79}. Also see \cite{Wang2001}.

\medskip

By the class field theory, idele class characters $\chi$ of finite order
correspond in a canonical way to characters of ${\rm Gal} (\bar{K} / K)$ of
finite order. By abuse of notation, we will still use the letter
$\chi$ to denote this Galois character. Moreover, there is a
canonically associated finite abelian extension $L / K$ such that the
kernel of $\chi$ as an idele class character is $N_{L / K} C_{L} \subset C_{L}$,
the norm from
$C_{L}$.

\medskip

Throughout this section, for any number field $L$, $d_{L}$ denotes
the discriminant of $L$, and $d_{L / K}$ denotes the relative
discriminant of $L$ over a subfield $K$.

\medskip

\subsection*{Preparations} \MED

\medskip

First we introduce two kernel functions used in the classical
analytic method, which was also used by Lagarias,
Odlyzko, Montgomery and K.\@ Murty (cf. \cite{L-O77},
\cite{L-M-O79}, \cite{KM94} and \cite{Serre81}). The use of these
two different kernel functions is related to the \emph{Explicit
formulas} of A.P.\@ Guinand (\cite{Guinand48}) and A.\@ Weil
(\cite{Weil52}).

\medskip

Let
\begin{align}
    k_{1} (s) &= k_{1} (s; x, y) =
        {\left( \frac{y^{s - 1} - x^{s - 1}}{s - 1} \right)}^{2}
    \notag \\
    k_{2} (s) &= k_{2} (s; x) = x^{s^{2} + s}
    \notag
\end{align}

Thus
\begin{align}
    k_{1} (1) &=
        {\left( \log \frac{y}{x} \right)}^{2}
    \notag \\
    k_{2} (1) &= x^{2}
    \notag
\end{align}

\medskip

For each smooth function $k (s)$, denote $\hat{k} (u)$ the inverse
Mellin transform, defined as
\[
    \hat{k} (u) = \frac{1}{2 \pi i} \int_{a - i \infty}^{a + i \infty}
    k(s) u^{- s} ds
\]
where $a$ is a sufficiently large number.

\medskip

Thus for $a > 1$, we have
\begin{align}
    \hat{k}_{1} (u) &= \hat{k}_{1} (u; x, y) =
        \begin{cases}
            0, &\text{if $u \ge y^{2}$ or $u \leq x^{2}$;} \\
            \frac{1}{u} \log \frac{y^{2}}{u} &\text{if $x y \leq u \leq y^{2}$;} \\
            \frac{1}{u} \log \frac{u}{x^{2}} &\text{if $x^{2} \leq u \leq x y$.}
        \end{cases}
    \notag \\
    \hat{k}_{2} (u) &= \hat{k}_{2} (u; x) =
    {(4 \pi \log x)}^{-\frac{1}{2}} \exp \left\{
        -\frac{{\left( \log \frac{u}{x} \right)}^{2}}{4 \log x}
        \right\}
    \notag
\end{align}

\medskip

Note that for each $j$ and $u$, $\hat{k}_{j} (u) \ge 0$,  and for
large $u$, $\hat{k}_{j} (u)$ is small.

\medskip

\begin{lemma} \label{T:301}
Assume that $L / K$ is cyclic.

\medskip

\textnormal{(1)} Let $\Sigma^{R}$ denote the summation over the prime
ideals of $K$ that ramify in $L$, then
\begin{align}
            &{\sum_{\fP}}^{R} \sum_{m \ge 1} \log (N \fP)
            \hat{k}_{1} (N \fP^{m})
            << \frac{\log \frac{y}{x}}{x^{2}} \log d_{L}
            \notag \\
            &{\sum_{\fP}}^{R} \sum_{m \ge 1, N \fP^{m} \leq x^{10}} \log (N \fP)
            \hat{k}_{2} (N \fP^{m})
            << {(\log x)}^{\frac{1}{2}} \log d_{L}
            \notag
\end{align}

\medskip

\textnormal{(2)} Let $\Sigma^{S}$ denote the summation over the prime
ideals of $K$ in $S$, then
\begin{align}
            &{\sum_{\fP}}^{S} \sum_{m \ge 1} \log (N \fP)
            \hat{k}_{1} (N \fP^{m})
            <<  \frac{\log \frac{y}{x}}{x^{2}} \log N_{S}
            \notag \\
            &{\sum_{\fP}}^{S} \sum_{m \ge 1, N \fP^{m} \leq x^{10}} \log (N \fP)
            \hat{k}_{2} (N \fP^{m})
            <<  {(\log x)}^{\frac{1}{2}} \log N_{S}
            \notag
\end{align}

\medskip

\textnormal{(3)} Let $\Sigma^{P}$ denote the summation over the pairs $(\fP, m)$
for which $N \fP^{m}$ is not a rational prime, then
\begin{align}
            &{\sum_{\fP, m}}^{P} \log (N \fP)
            \hat{k}_{1} (N \fP^{m})
            << n_{K} \frac{\left(\log \frac{y}{x} \right) (\log y)}{x (\log x)}
            \notag \\
            &{\sum_{\fP, m}}^{P} \log (N \fP)
            \hat{k}_{2} (N \fP^{m})
            << n_{K} x^{7/4}
            \notag
\end{align}

\medskip

\textnormal{(4)}
\[
            {\sum_{\fP}} \sum_{m \ge 1, N \fP^{m} > x^{3 + \delta}} \log (N \fP)
            \hat{k}_{2} (N \fP^{m})
            << n_{K} x^{2 - \frac{\delta^{2}}{4}} (\log x)
\]
where $\delta$ is any positive number.
\end{lemma}

\medskip

\emph{Proof.} (Also see \cite{L-O77} and \cite{L-M-O79})

\medskip

(1-a, cf.\ Lemma 3.1 of \cite{L-M-O79}):
\begin{align}
\sum_{\fP}^{R} & \sum_{m = 1}^{\infty} \log (N \fP) \hat{k}_{1} (N \fP^{m}) \notag\\
&<< \sum_{\fP}^{R} \log (N \fP) \log (y / x) \sum_{m = 1}^{\infty} {(N \fP)}^{-m} \notag\\
&<< \frac{\log (y / x)}{x^{2}} \sum_{\fP}^{R} \log (N \fP) \notag\\
&<< \frac{\log (y / x)}{x^{2}} \log (d_{L})
\notag
\end{align}

\medskip

Here we use two facts: (i) $\hat{k}_{1} (u) << u^{-1} \log (y / x)$ and $\hat{k}_{1} = 0$
if $u < x^{2}$. (ii) The conductor-discriminant formula
$\sum_{\fP}^{R} \log (N \fP) << \log (d_{L})$.

\medskip

(1-b, cf.\ Lemma 3.1 of \cite{L-M-O79}):
\begin{align}
\sum_{\fP}^{R} & \sum_{m \ge 1, N \fP^{m} \leq x^{10}} \log (N \fP) \hat{k}_{2} (N \fP^{m})
\notag \\
&<< \sum_{\fP}^{R} \log (N \fP) \sum_{m \ge 1, N \fP^{m} \leq x^{10}} {(\log x)}^{-1/2}
\notag\\
&<< \sum_{\fP}^{R} \log (N \fP) {(\log x)}^{1/2} \notag\\
&<< {(\log x)}^{1/2} \log (d_{L}) \notag
\end{align}

\medskip

Here we use two trivial estimates $\hat{k}_{2} (u) << {(\log x)}^{-1/2}$
and $\# \{ N \fP^{m} \leq x^{10}\} << \log x$.

\medskip

(2) the proof is
almost the same as (1) except that we need to estimate\linebreak $\sum^{S}
\log (N \fP)$ instead of $\sum^{R} \log (N \fP)$.

\medskip

(3-a, cf.\ see Lemma 3.2 in \cite{L-M-O79}) Use the fact that the number of pairs
$(\fP, m)$ such that $N \fP^{m} = q$ is at most $n_{K}$, and we have
\begin{align}
    {\sum_{\fP, m}}^{P} \log (N \fP) \hat{k}_{1} (N \fP^{m})
        &<< n_{K} (\log \frac{y}{x}) \sum_{x^{2} \leq p^{h} \leq y^{2}, h \ge 2}
        p^{-h} \log p^{h}
        \notag \\
        &<< n_{K} (\log \frac{y}{x}) (\log y) \sum_{n = p^{a}, a \ge 2, n \ge x^{2}}
        n^{-1}
        \notag \\
        &<< n_{K} (\log \frac{y}{x}) (\log y) \frac{1}{x \log x}
        \notag
\end{align}
where the last bound uses the prime number theorem.

\medskip

(3-b, cf.\ see Lemma 3.2 in \cite{L-M-O79}) let $S (u)$ denote the number of prime
power integers $p^{h} (h \ge 2)$ in the interval $[1, u]$. It is easy to see
that $S (u) << u^{1/2}$ since
$S (u) \leq u^{1/2} + u^{1/3} + u^{1/5} + \cdots << u^{1/2} + u^{1/3} \log u$. Thus
\begin{align}
    {\sum_{\fP}}^{P} \log (N \fP) \hat{k}_{2} (N \fP^{m})
    &<< n_{K} \sum_{p, h \ge 2} \log (p^{h}) \hat{k}_{2} (p^{h})
    \notag \\
    &<< n_{K} \int_{3}^{\infty} (\log u)\,\hat{k}_{2} (u)\,dS (u)
    \notag \\
    &< -n_{K} \int_{3}^{\infty} S (u) d\, (\log (u) \hat{k}_{2} (u))
    \notag \\
    &<< n_{K} \int_{3}^{\infty} u^{1/2} \log (u) (-\hat{k}'_{2} (u)) d u
    \notag \\
    &<< n_{K} \int_{3}^{\infty} u^{1/2} \log (u) \hat{k}_{2} (u)
    \frac{\log (u / x)}{2 \log x} \frac{x}{u} d u
    \notag \\
    &<< n_{K} \frac{x}{{(\log x)}^{3/2}} \int_{3}^{\infty} u^{-1/2} {(\log u)}^{2}
    \exp \left( - \frac{{(\log (u / x))}^{2}}{4 \log x} \right) d u
    \notag
\end{align}
\begin{align}
    &<< n_{K} \frac{x^{3/4}}{{(\log x)}^{3/2}} \int_{3}^{\infty} {(\log u)}^{2}
    \exp \left( - \frac{{(\log u)}^{2}}{4 \log x} \right) d u
    \notag \\
    &<< n_{K} \frac{x^{3/4}}{{(\log x)}^{3/2}} \int_{-\infty}^{\infty} t^{2}
    \exp \left( - \frac{t^{2}}{4 \log x} + t \right) d t
    \notag \\
    &<< n_{K} x^{3/4} \int_{-\infty}^{\infty} t^{2}
    \exp \left( -t^{2} + t \sqrt{\log x} \right) d t
    \notag \\
    &<< n_{K} x^{3/4} \int_{-\infty}^{\infty} {(t + \frac{\sqrt{\log x}}{2})}^{2}
    \exp (-t^{2}) x^{1/4} d t
    \notag \\
    &<< n_{K} x \log (x) < n_{K} x^{7/4} \notag
\end{align}

\medskip

For (4),
\begin{align}
    \sum_{N \fP^{m} > x^{3 + \delta}} & \log (N \fP) \hat{k}_{2} (N \fP^{m})
    << n_{K} \sum_{q > x^{3 + \delta}} (\log q) \hat{k}_{2} (q)
    \notag \\
    &<< n_{K} \int_{x^{3 + \delta}}^{+\infty} (\log u) \hat{k}_{2} (u) du
    \notag \\
    &<< n_{K} {(\log x)}^{-\frac{1}{2}} \int_{(3 + \delta) \log x}^{+\infty}
    t \exp \left( -\frac{{(t - \log x)}^{2}}{4 \log x} \right) e^{t} dt
    \notag \\
    &<< n_{K} x^{2} {(\log x)}^{-\frac{1}{2}} \int_{(3 + \delta) \log x}^{+\infty}
    t \exp \left( -\frac{{(t - 3 \log x)}^{2}}{4 \log x} \right) dt
    \notag
\end{align}
\begin{align}
    &<< n_{K} x^{2} {(\log x)}^{-\frac{1}{2}}
    \Bigl\{
        \int_{\delta \log x}^{+\infty} 3 \log x
        \exp \left( -\frac{t^{2}}{4 \log x} \right) dt
    \notag \\
    & \qquad + (\log x) \int_{\delta \sqrt{\log x}}^{+\infty}
        t \exp \left( -\frac{t^{2}}{4} \right) dt
    \Bigr\}
    \notag \\
    &<< n_{K} x^{2} {(\log x)}^{-\frac{1}{2}}
    \Bigl\{
        3 {(\log x)}^{3/2} \int_{\delta \sqrt{\log x}}^{+\infty}
        \exp \left( - \frac{ t^{2}}{4} \right) dt
    \notag \\
    & \qquad + (\log x) \exp \left( - \frac{\delta^{2}}{4} \log x \right)
    \Bigr\}
    \notag \\
    &<< n_{k} x^{2 - \frac{\delta^{2}}{4}} (\log x)
    \notag
\end{align}
where we use the following: $\hat{k}_{2} (u) << \hat{k}_{2} (x)$
if $|x - u| < 1$, and a well known estimate $\int_{T}^{+\infty} e^{-
t^{2}/4} dt << e^{- T^{2}/4}$.

\medskip

\qedsymbol \medskip

\begin{lemma} \label{T:302}
Let $\chi$ be a global character of $C_{K}$.

\medskip

\textnormal{(1)} If $N (t) = N_{L} (t)$ denotes the number of zeros $\rho = \beta + i \gamma$,
of $\zeta_{L} (s)$ with $0 < \beta < 1$ and $|\gamma - t| \leq 1$,
then we have
\[
            N (t) << \log d_{L} + n_{L} \log (|t| + 2)
\]

\medskip

\textnormal{(2)} If $n (r; s) = n_{L} (r; s)$ denotes the number of zeros $\rho$,
of $\zeta_{L} (s)$ with $|\rho - s| \leq r$, then we have
\[
            n (r; s) << 1 + r (\log d_{L} + n_{L} \log (|s| + 2))
\]

\medskip

\textnormal{(3)} If $N_{\chi} (t)$ denotes the number of zeros $\rho = \beta + i \gamma$,
of $L(s, \chi, K)$ with $0 < \beta < 1$ and $|\gamma - t| \leq 1$, then we have
\[
            N_{\chi} (t) << \log A (\chi) + n_{K} \log (|t| + 2)
\]

\medskip

\textnormal{(4)} If $n_{\chi} (r; s)$ denotes the number of zeros $\rho$
of $L(s, \chi, K)$ with $|\rho - s| \leq r$, then we have
\[
            n_{\chi} (r; s) << 1 + r (\log A (\chi) + n_{K} \log (|s| + 2))
\]
\end{lemma}

\emph{Proof. } These are standard results and the proof can be found
in a lot of literatures. For examples, see \cite{Neukirch91},
\cite{L-O77} and Lemma 2.2 in
\cite{L-M-O79}.

\medskip

\qed \medskip

The following is the Siegel zero free region result for $L (s, \chi)$.

\begin{lemma} \label{T:303}
Let $\chi$ be a global character of $K$. There is a positive, absolute,
effectively computable constant $c_{2}$ such that

\medskip

\textnormal{(1)} $L(s, \chi)$ has no zero $\rho = \beta + i \gamma$
        in the region
        \begin{align}
            \beta &\ge 1 - {c_{2}}^{-1}
             {(\log A(\chi) + n_{K} \log (|\gamma| + 2))}^{-1}
            \notag \\
            \gamma &\ge {(1 + c_{2} \log A(\chi))}^{-1}
            \notag
        \end{align}
    where $\gamma \ne 0$.

\medskip

\textnormal{(2)} $L(s, \chi)$ has at most one zero in the region
        \begin{align}
            \beta &\ge 1 - {(c_{2} \log A(\chi))}^{-1}
            \notag \\
            \gamma &\leq {c_{2} \log A(\chi))}^{-1}
            \notag
        \end{align}

\medskip

If such a zero exists, it must be simple and real, and $\chi$ must be trivial or quadratic.
\end{lemma}

\medskip

\emph{Proof. } See \cite{L-O77}, or \cite{Lang70} and
\cite{Neukirch91}. \qed

\medskip

Before finishing this part, we quote the Deuring--Heilbronn
phenomenon here, a discussion of which can be found in Section 5 in
\cite{L-M-O79}.

\medskip

\begin{lemma} \textnormal{\bf[Deuring--Heilbronn Phenomenon]} \label{T:304}

    There are positive, absolute,   effectively computable constants $c_{7}$
    and $c_{8}$ such that if $\zeta_{L} (s)$ has a real zero $\beta_{0}$,
    then $\zeta_{L} (\sigma + it) \ne 0$ for
    \[
        \sigma \ge 1 - c_{8} \cdot \frac{
            \log \left(
                \frac{c_{7}}{(1 - \beta_{0}) \log (d_{L} \tau^{n_{L}})}
            \right)
        }{\log (d_{L} \tau^{n_{L}})}
    \]
    where $\tau = |t| + 2$ with the single exception $\sigma + it = \beta_{0}$.
\end{lemma}

\medskip \qed \medskip

\begin{corollary} \label{T:305}
    There is a positive, absolute, effectively computable constant $c_{10}$
    such that any real zero $\beta_{0}$ of $\zeta_{L} (s)$ satisfies
    \[
        1 - \beta_{0} \ge {d_{L}}^{-c_{10}}
    \]
\end{corollary}

\medskip \emph{Proof.} See Corollary 5.2 in \cite{L-M-O79}. \qed

\medskip

\subsection*{Standard Model}  \MED

In this part, we will recall the main model of
\cite{L-M-O79} for our method here. We have included the relevant
details for the convenience of the readers.

\medskip

We need to consider the Artin $L$--series $L(s, \phi, L/K)$ (cf.
\cite{Lang70},
\cite{Neukirch91}, \cite{L-O77}, \cite{L-M-O79}) where $\phi$ is the
character of an irreducible representation of $G = {\rm Gal} (L/K)$. We have
\[
    -\frac{L'}{L} (s, \phi, L/K) =
    \sum_{\fP} \sum_{m \ge 1} \Phi_{K} (\fP^{m}) \log (N \fP) {(N \fP)}^{-m s}
\]
where
\[
    \Phi_{K} (\fP^{m}) = \frac{1}{e_{\fP}(L/K)} \sum_{\alpha \in I_{\fP}(L/K)}
        \phi (\tau^{m} \alpha)
\]
where $\tau = \left( \frac{L/K}{\fP} \right)$ is one representative
of the Frobenius element corresponding to $\fP$, $I_{\fP} =
I_{\fP}(L/K)$ is the inertial subgroup of the decomposition group
$G_{\fP} = {\rm Gal} (L_{\fQ}/K_{\fP})$ and $e_{\fP}(L/K) = |I_{\fP}|$ is the
ramification index of $\fQ$ over $\fP$.

\medskip

If $\fP$ is unramified in $L$ then $\Phi_{K} (\fP^{m}) = \phi
(\alpha^{m})$. If $L / K$ is abelian, then all irreducible $\phi$ are
characters (and hence by the class field theory, correspond
to Hecke characters).

\medskip

\begin{lemma} \label{T:306}
    Let $C$ be a conjugacy class of $G$ and $g$ a representative of $C$,
    $H = <g>$ and $E = L^{H}$ the fixed field of $g$. Then we have

    \textnormal{(1)}
    \[
        F_{C} (s) := - \frac{|C|}{|G|}
        \sum_{\phi \text{ irreducible}} \bar{\phi} (g)
        \frac{L'}{L} (s, \phi, L/K)
        =  - \frac{|C|}{|G|}
        \sum_{\chi \in \hat{G} (L/E)} \bar{\phi} (g)
        \frac{L'}{L} (s, \chi, E)
    \]
    where $\hat{G} (L/E)$ denotes the group of characters
    of $G(L/E)$,
    and
    \textnormal{(2)}
    \[
        F_{C} (s) = \sum_{\fP} \sum_{m \ge 1} \theta (\fP^{m})
        \log (N \fP) {(N \fP)}^{-s},
    \]
    where
    \[
        \theta (\fP^{m}) = \begin{cases}
            1 &\text{if ${\left( \frac{L/K}{\fP} \right)}^{m} = C$} \\
            0 &\text{if ${\left( \frac{L/K}{\fP} \right)}^{m} \ne C$}
        \end{cases}
    \]
    and $|\theta (\fP^{m})| \leq 1$ if $\fP$ ramifies in $L$.
\end{lemma}

\medskip \emph{Proof. } This is an exercise
of representations theory. See Section 5, \cite{L-O77}. \qed

\medskip

The previous lemma allows us to reduce the density problem to the
case of a cyclic extension, for which we can use just the abelian
$L$--series of Hecke.

\medskip

The following lemma (cf. \cite{Lang70}, \cite{L-O77}) describes a
functional equation that \linebreak
$L(s, \chi, E)$ satisfies.

\medskip

\begin{lemma} \label{T:307}
    Let $L(s, \chi) = L(s, \chi, E)$ be the $L$--series associated to $\chi \in \hat{G} (L/E)$.
    \[
        A (\chi) = d_{E} N_{E/\bQ} (\fF_{0} (\chi))
    \]
    where $\fF_{0} (\chi)$ denotes the finite conductor of $\chi$.
    \[
        \delta (\chi) = \begin{cases}
            1 &\text{if $\chi$ is principal;} \\
            0 &\text{otherwise.}
        \end{cases}
    \]

    There are nonnegative integers $a = a (\chi)$ and $b = b (\chi)$
    such that
    \[
        a (\chi) + b (\chi) = n_{E}
    \]

    Set
    \[
        \gamma_{\chi} (s) = {\left\{ \pi^{- \frac{s + 1}{2}}
        \Gamma (\frac{s + 1}{2}) \right\} }^{b}
        {\left\{ \pi^{- \frac{s}{2}}
        \Gamma (\frac{s}{2}) \right\} }^{a}
    \]
    and
    \[
        \Lambda (s, \chi) = {(s (s - 1))}^{\delta (\chi)}
        {A (\chi)}^{s/2} \gamma_{\chi} (s) L(s, \chi)
    \]

    Then $\Lambda (s, \chi)$ satisfies the functional equation
    \[
        \Lambda (s, \chi) = W (\chi) \Lambda (1 - s, \bar{\chi})
    \]
    where $W (\chi)$ is a certain constant of absolute $1$.

    Furthermore, $\Lambda (s, \chi)$ is entire of order $1$ and does not vanish at
    $s = 0$.
\end{lemma}

\qed

\medskip

Let
\[
    J_{j} (\chi) \stackrel{\triangle}{=} - \frac{1}{2 \pi i} \int_{2 - i \infty}^{2 + i \infty}
    \frac{L'}{L} (s, \chi) k_{j} (s) ds
\]
and
\[
    I_{j} = - \frac{1}{2 \pi i} \int_{2 - i \infty}^{2 + i \infty}
    F_{C} (s) k_{j} (s) ds
\]
where $F_{C} (s)$ is defined in Lemma \ref{T:306}.

\medskip

By this lemma, we have
\begin{equation}  \label{E:3-01} \tag{3-1}
    I_{j} = \frac{|C|}{|G|} \sum_{\chi \in \hat{G} (L/E)} \bar{\chi} (g) J_{j} (\chi)
\end{equation}
where $g$ is a representative of $C$.

\medskip

We have two ways to express $I_{j}$. One way is using the inverse
Mellin transform and the other is using the residue theorem.

\medskip

By the inverse Mellin transform, and we have
\[
    J_{j} (\chi) = \sum_{\fP} \sum_{m \ge 1} \chi (\fP^{m}) \log (N \fP)
        \hat{k}_{j} (N \fP^{m})
\]
since
\[
    - \frac{L'}{L} (s, \chi, E) = \sum_{\fP} \sum_{m \ge 1} \chi (\fP^{m}) \log (N \fP)
         {(N \fP)}^{- m s}
\]

\medskip

Also, by Lemma ~\ref{T:306},
\begin{align}
    I_{j} &= \frac{|C|}{|G|} \sum_{\chi \in \hat{G} (L/E)} \bar{\chi} (g)
        \sum_{\fP} \sum_{m \ge 1} \chi (\fP^{m}) \log (N \fP)
        \hat{k}_{j} (N \fP^{m})
    \notag \\
    &= \sum_{\fP} \sum_{m \ge 1} \theta (\fP^{m}) \log (N \fP)
        \hat{k}_{j} (N \fP^{m})
    \notag
\end{align}

\medskip

\begin{lemma} \label{T:308} \MED

    \noindent \textnormal{(1)}
    \begin{align}
        J_{j} (\chi) &= \delta (\chi) k_{j} (1) - \sum_{\rho} k_{j} (\rho)
        + O (n_{E} k_{j} (0))
        \notag \\
        &+ O (k_{j} (- \frac{1}{2}) (\log A (\chi) + n_{E}))
        \notag
    \end{align}
    where the sum runs over all the nontrivial zeros of
     $L(s, \chi, E)$,
    and all the implied constants are absolute and effectively computable.

    \medskip

    \noindent \textnormal{(2)}
    \begin{align}
        \frac{|G|}{|C|} I_{j}
        &\ge  k_{j} (1) - \sum_{\rho} k_{j} (\rho)
        \notag \\
        &- c_{6} \{ n_{L} k_{j} (0) + k_{j} (- \frac{1}{2}) \log d_{L} \}
        \notag
    \end{align}
    where the sum runs over all the nontrivial zeros of $\zeta_{L} (s)$
    and $c_{6}$ is positive, absolute and effectively computable.
\end{lemma}

\medskip

For the proof we need the following

\begin{proposition} \textnormal{(\bf
     the Conductor--Discriminant Formula)}  \label{T:309}
    \[
        \prod_{\chi \in \hat{G} (L/E)} A (\chi) = d_{L}
    \]
\end{proposition}
\qed

\medskip

For a proof, see \cite{L-O77}, \cite{Odlyzko77}.

\medskip

\emph{Proof of Lemma ~\ref{T:308}.}

For (1), see \cite{L-O77}. The basic idea is to consider the
following integral
\begin{align}
    J_{j} (\chi, T) &\stackrel{\triangle}{=} - \frac{1}{2 \pi i}
    \int_{\partial B (T)} \frac{L'}{L} (s, \chi, E) k_{j} (s) ds
    \notag \\
    &= \delta (\chi) k_{j} (1) - a_{\chi} k_{j} (0) - \sum_{|\gamma| < T} k_{j} (\rho)
    \notag
\end{align}
where the sum runs over all the zeros $\rho = \beta + i \gamma$ of
$L(s, \chi, E)$ within the rectangle $B(T)$: $[- \frac{1}{2}, 2]
\times [- T, T]$. Estimate the integral on each line segment and let
$T$ go to the infinity as in \cite{L-O77}. In fact, on the line
segment from $-\frac{1}{2} + iT$ to $-\frac{1}{2} - iT$,
\[
    \left| \frac{L'}{L} (s, \chi, E) \right|
    << \log A (\chi) + n_{E} (\log (|s| + 2))
\]
(see Lemma 6.2, \cite{L-O77}). Thus,
\[
    \left|
        \frac{1}{2 \pi i} \int_{-\frac{1}{2} - iT}^{-\frac{1}{2} + iT}
        \frac{L'}{L} (s, \chi, E) k_{j} (s)ds
    \right|
    << k_{j} (-\frac{1}{2}) \{ \log A(\chi) + n_{E} \}
\]
as
\begin{align}
    k_{1} (- \frac{1}{2} + i t) &<< k_{1} (- \frac{1}{2}) \frac{1}{1 + t^{2}}
    \quad \text{if $y >> x$}
    \notag \\
    k_{2} (- \frac{1}{2} + i t) &<< k_{1} (- \frac{1}{2}) \exp (-t^{2} \log x)
    \quad \text{if $x >> 1$}
    \notag
\end{align}

\medskip

To estimate the integral $I_{\pm} (T)$ on the horizontal line
segments from $2 \pm iT$ to $-\frac{1}{2} \pm iT$, one uses the
method of Landau (Section 6 of \cite{L-O77}, Section 3 of
\cite{L-M-O79}, \cite{Landau27}), obtaining the estimate
\[
    I_{\pm} (T) << k_{j} (iT) (\log A (\chi) + n_{E} \log T).
\]
Note that $T \to \infty$, $I_{\pm} (T) \to 0$.

\medskip

Combining these estimates with Proposition ~\ref{T:309}, we obtain
(1).

\medskip

Now (2) is easy to get from (1) since
\[
    I_{j} = \frac{|C|}{|G|} \sum_{\chi \in \hat{G} (L/E)} J_{j} (\chi)
\]
and
\[
    \zeta_{L} (s) = \prod_{\chi \in \hat{G} (L/E)} L(s, \chi, E)
\]
and we can use Proposition ~\ref{T:309}.

\medskip

\qed

\medskip

Now we are ready to explain how we plan to use the standard model
for our purposes.

From the rest of this chapter, assume that $y >> x$ if we apply the
first kernel function $k_{1} (s)$ and $x >> 1$ if we apply the
second one $k_{2} (s)$. Let $n = n_{L} / n_{K}$ which is not less
than $|G|/|C|$.

Thus, by Lemma ~\ref{T:308} (2), we have
\[
    I_{j} \ge \frac{1}{n} (k_{j} (1) - \sum_{\rho} |k_{j} (\rho)|)
    - c_{6} \left\{ n_{K} k_{j} (0) + k_{j} (-\frac{1}{2})
    \left( \frac{1}{n} \log d_{L} \right) \right\}
\]

Note that
\begin{align}
    k_{1} (0) &= {\left( \frac{x^{-1} - y^{-1}}{-1} \right)}^{2} << x^{-2}
    \notag \\
    k_{1} (-\frac{1}{2}) &= {\left( \frac{x^{-\frac{3}{2}} - y^{-\frac{3}{2}}}{-\frac{3}{2}}
        \right)}^{2} << x^{-3}
    \notag \\
    k_{2} (0) &= 1
    \notag \\
    k_{2} (-\frac{1}{2}) &= x^{-\frac{1}{4}}
    \notag
\end{align}
Thus, the $c_{6} \{\,\*\,\}$ term is bounded by some multiple of
\[
    T_{j} = \begin{cases}
        \frac{x^{-2}}{n} \log d_{L}, &\text{if $j = 1$;} \\
        \frac{1}{n} \log d_{L}, &\text{if $j = 2$.}
    \end{cases}
\]

\medskip

Furthermore, we have
\[
    I_{j} = \sum_{\fP} \sum_{m \ge 1} \theta (\fP^{m}) \log (N \fP)
        \hat{k}_{j} (N \fP^{m})
\]

\medskip

Thus,
\begin{align}
    I_{1} &= S_{1, 1} + S_{1, 2} + S_{1, 3} + \tilde{I}_{1}
    \notag \\
    I_{2} &= S_{2, 1} + S_{2, 2} + S_{2, 3} + S_{2, 4}
    + \tilde{I}_{2}
    \notag
\end{align}
where the symbols mean the following:

\medskip

$\tilde{I}_{j}$ denotes the sum over the primes outside $S$,
unramifying in $L$, of degree $1$ over $K$ and the Artin symbol of
$\fP$ under $L/K$ being $C$ such that $N \fP \leq y^{2}$ or $x^{3 +
\delta}$ when $j = 1$ or $2$ respectively.

\medskip

$S_{1, 1}$ denotes the sum over $(\fP, m)$ with $\fP$ ramifying in
$L$. $S_{2, 1}$ denotes the sum over $(\fP, m)$ with $\fP$ ramifying
in $L$ and $N \fP^{m} \leq x^{10}$.

\medskip

$S_{1, 2}$ denotes the sum over $(\fP, m)$ with $\fP$ in $S$. $S_{2,
2}$ denotes the sum over $(\fP, m)$ with $\fP$ in $S$ and $N \fP^{m}
\leq x^{10}$.

\medskip

$S_{j, 3}$ denotes the sum over $(\fP, m)$ with $N \fP^{m}$ not a
rational prime.

\medskip

$S_{2, 4}$ denotes the sum over $(\fP, m)$ with $N \fP^{m} > x^{3 +
\delta}$.

\medskip

Applying Lemma ~\ref{T:301}, we have
\begin{align}
    S_{1, 1} &<< \frac{1}{n} \frac{\log (y/x)}{x^{2}} \log d_{L}
    \notag \\
    S_{2, 1} &<< \frac{1}{n} {(\log x)}^{\frac{1}{2}} \log d_{L}
    \notag \\
    S_{1, 2} &<< \frac{\log (y/x)}{x^{2}} \log N_{S}
    \notag \\
    S_{2, 2} &<< {(\log x)}^{\frac{1}{2}} \log N_{S}
    \notag \\
    S_{1, 3} &<< n_{K} \frac{(\log (y/x))(\log y)}{x (\log x)}
    \notag \\
    S_{2, 3} &<< n_{K} x^{7/4}
    \notag \\
    S_{2, 4} &<< n_{K} x^{2 - \frac{\delta^{2}}{4}} \log x
    \notag
\end{align}

Then the main idea of this model is the following: Pick $x$, $y$
appropriately. If we assume that for any $\fP$ unramifying in $L$,
of degree $1$ over $K$ and the Artin symbol of $\fP$ under $L/K$
being $C$ such that either $N \fP > y^{2}$ or $x^{3 + \delta}$ when
$j = 1$ or $2$ respectively, or $\fP \in S$ or $\fP$ ramifies in
$L$, then $\tilde{I}_{j} = 0$ and
\[
    \frac{1}{n} \left( k_{j} (1) - \sum_{\rho} |k_{j} (\rho)| \right)
    \leq c'_{6} T_{j} + \sum_{v} S_{j, v}
\]
However, if the left--hand side dominates over $c'_{6} T_{j}$ and
$S_{j, v}$ by a sufficiently large constant factor, then one gets a
contradiction.

\medskip

So the key component of this model is to find a better lower bound
for
\[
    k_{j} (1) - \sum_{\rho} |k_{j} (\rho)|
\]

\medskip

\subsection*{Final Estimations}  \MED

In this part we will prove Theorem ~\ref{TM:A}. Let $P_{1}
(C, S)$ be the set of primes of $K$ satisfying (1) to (3) in Theorem
~\ref{TM:A}.

\medskip

From last section we've already seen that the quality of the
effective bound depends on the lower bound of $k_{j} (1) -
\sum_{\rho} |k_{j} (\rho)|$. However, the possible exceptional zero
$\beta_{0}$ will cause difficulty. In general, one will be forced to
use the Deuring--Heilbronn. Fortunately, there is nothing new here
compared with the classical case where $S = \emptyset$.

\medskip

To simplify our notation, we define $\beta_{0}$ to be the
exceptional zero of $\zeta_{L} (s)$ if it exists, and $\beta_{0} = 1
- {(c_{2} \log d_{L})}^{-1}$ otherwise, where $c_{2}$ is the
constant defined in Lemma ~\ref{T:303}, so that $\zeta_{L} (s)$ has
at most one zero in the interval $(\,1 - {(c_{2} \log d_{L})}^{-1},
1\,)$.

\medskip

In either case,
\[
    k_{j} (1) - \sum_{\rho} |k_{j} (\rho)| \ge
    k_{j} (1) - k_{j} (\beta_{0}) - \sum_{\rho \ne \beta_{0}} |k_{j} (\rho)|
\]

\medskip

By using the mean value theorem, we have
\begin{align}
    k_{1} (1) - k_{1} (\beta_{0}) &= {\left( \log \frac{y}{x} \right)}^{2}
        - {\left( \frac{y^{\beta_{0} - 1} - x^{\beta_{0} - 1}}{\beta_{0} - 1} \right)}^{2}
    \notag \\
    &\ge \frac{1}{10} {\left( \log \frac{y}{x} \right)}^{2}
    \min \left\{ \,1, (1 - \beta_{0}) \log \left( \frac{y}{x} \right) \,\right\}
    \notag \\
    k_{2} (1) - k_{2} (\beta_{0}) &= x^{2} - x^{\beta_{0} + {\beta_{0}}^{2}}
    \ge \frac{x^{2}}{10} \min \left\{ \,1, (1 - \beta_{0}) \log (x) \,\right\}
    \notag
\end{align}

\medskip

First suppose
\[
    1 - \beta_{0} \ge {c_{7}}^{2} {(\log d_{L} 3^{n_{L} N_{S}^{n}})}^{-2}
\]
where $c_{7}$ is the constant defined in Lemma ~\ref{T:304}. In this
case, we use the kernel $k_{1} (s)$. (Recall that $n = n_{L} / n_{K}$.)

\medskip

The contribution of the zeros $\rho$ of $\zeta_{L} (s)$ with $|\rho
- 1| \ge 1$ is bounded by
\[
    \sum_{|\rho - 1| \ge 1} |k_{1} (\rho)| \leq \int_{1}^{\infty}
    \frac{2}{t^{2}} dn(t; 1) << \log d_{L}
\]
where $n(t; 1)$ is the number of the nontrivial zeros of $\zeta_{L}$
with $|\rho - 1| \ge t$ (See Lemma ~\ref{T:302}).

\medskip

Next, assume that $|\rho - 1| \leq 1$ for a nontrivial zero $\rho =
\beta + i \gamma \ne \beta_{0}$ of $\zeta_{L}$.

\medskip

If $\beta_{0}$ as an exceptional zero exists with
\[
    1 - \beta_{0} \leq \frac{1}{18} c_{2} {c_{7}}^{2} {(\log (d_{L} N_{S}^{n}))}^{-1},
\]
then since $d_{L} \ge 3^{n_{L}/2}$ for $n_{L} \ge 2$, we have
\[
    \frac{c_{7}}{(1 - \beta_{0}) \log (d_{L} 3^{n^{L}} N_{S}^{n})}
    \ge {\left\{
        (\frac{1}{2} c_{2}) (1 - \beta_{0}) \log (d_{L} N_{S}^{n})
    \right\} }^{-\frac{1}{2}},
\]
and therefore by the Deuring--Heilbronn (Lemma ~\ref{T:305}, note that when replace
$d_{L}$ by any $Q > d_{L}$, we still have this statement),
\[
    \beta \leq 1 - c_{8} \frac{ \log \left\{
        \frac{c_{7}}{(1 - \beta_{0}) \log (d_{L} 3^{n_{L}} N_{S}^{n})}
    \right\} }{\log (d_{L} 3^{n_{L}} N_{S}^{n})}
    \leq
    1 - c_{11}
    \frac{ \log {\left\{
    (\frac{1}{2} c_{2}) (1 - \beta_{0}) \log (d_{L} N_{S}^{n})
    \right\}}^{-1} }{\log (d_{L} N_{S}^{n})}
\]

\medskip

On the other hand, if
\[
    1 - \beta_{0} \ge \frac{1}{18} c_{2} {c_{7}}^{2} {(\log (d_{L} N_{S}^{n}))}^{-1}
\]
then by the zero--free region given by Lemma ~\ref{T:303},
\[
    \beta \leq 1 - {(3 c_{2} \log (d_{L}) N_{S}^{n})}^{-1}
\]

\medskip

Hence we have
\begin{equation} \label{*} \tag{*}
    \beta \leq 1 - c_{12}
    \frac{ \log {\left\{
    (\frac{1}{2} c_{2}) (1 - \beta_{0}) \log (d_{L} N_{S}^{n})
    \right\}}^{-1} }{\log (d_{L} N_{S}^{n})}
\end{equation}
for some $0 < c_{12} < c_{11}$.

\medskip

Thus \eqref{*} holds for all the cases.

\medskip

Let
\[
    B = c_{12}
    \frac{ \log {\left\{
    (\frac{1}{2} c_{2}) (1 - \beta_{0}) \log (d_{L} N_{S}^{S})
    \right\}}^{-1} }{\log (d_{L} N_{S}^{n})}
\]
From \eqref{*}, we have
\[
    |k_{1} (\rho)| << x^{2 (\beta - 1)} {|\rho - 1|}^{- 2}
    << x^{- 2 B} {|\rho - 1|}^{-2}
\]
Thus, by Lemma ~\ref{T:301},
\begin{align}
    \sum_{|\rho - 1| < 1, \rho \ne \beta_{0}} |k_{1} (\rho)|
    &\leq x^{- 2 B} \int_{B}^{1} \frac{1}{t^{2}} d n(t; 1)
    \notag \\
    &<< x^{- 2 B} (B^{-2} + B^{-1} \log d_{L})
    \notag \\
    &<< x^{- 2 B} B^{-1} \log d_{L}
    \notag
\end{align}
As $B >> {(\log (d_{L} N_{S}^{n}))}^{-1}$, using the expression of $B$, we have
\[
    \sum_{|\rho - 1| < 1, \rho \ne \beta_{0}} |k_{1} (\rho)|
    << {(\log (d_{L} N_{S}^{n}))}^{2} {
    \left\{
    (\frac{1}{2} c_{2}) (1 - \beta_{0}) \log (d_{L} N_{S}^{n})
    \right\} }^{2 c_{12} \frac{\log x}{\log (d_{L} N_{S}^{n})}}
\]
Thus we have shown that
\begin{align}
    \label{4A-1} \tag{4A-1}
    k_{1} (1)- \sum_{\rho} |k_{1} (\rho)|
    & \ge \frac{1}{10} {\left( \log \frac{y}{x} \right)}^{2}
    \min \left\{ \,1, (1 - \beta_{0}) \log \frac{y}{x} \,\right\}
    \\
    & - c_{13} \log d_{L}
    \notag \\
    & - c_{14} {(\log d_{L})}^{2} {
    \left\{
    (\frac{1}{2} c_{2}) (1 - \beta_{0}) \log d_{L}
    \right\} }^{2 c_{12} \frac{\log x}{\log d_{L}}}
    \notag
\end{align}
for some positive constants $c_{13}$ and $c_{14}$.

\medskip

We now complete the proof of Theorem ~\ref{TM:A} in the case
\[
    1 - \beta_{0} \ge {c_{7}}^{2} {(\log d_{L} 3^{n_{L}})}^{-2}.
\]

\medskip

Assume that for any $\fP$ in $P_{1} (C, S)$, $N \fP > y^{2}$. Then
\begin{align}
    0 &= \tilde{I}_{1} = \sum_{\fP \in P_{1}(C, S)}
    (\log N \fP) \hat{k}_{1} (N \fP)
    \notag \\
    & \ge \frac{1}{10 n} {\left( \log \frac{y}{x} \right)}^{2}
    \min \left\{ \,1, (1 - \beta_{0}) \log \frac{y}{x} \,\right\}
    \notag \\
    & - c_{13} \frac{1}{n} \log d_{L}
    \notag \\
    & - c_{14} \frac{1}{n} {(\log (d_{L} N_{S}^{n}))}^{2} {
    \left\{
    (\frac{1}{2} c_{2}) (1 - \beta_{0}) \log (d_{L} N_{S}^{n})
    \right\} }^{2 c_{12} \frac{\log x}{\log (d_{L} N_{S}^{n})}}
    \notag \\
    & - c_{15, 1} \left\{
    \frac{1}{n} \frac{1}{x^{2}} \log \left( \frac{y}{x} \right) \log d_{L}
    \right\}
    \notag \\
    & - c_{15, 2} \left\{
    \frac{1}{x^{2}} \log \left( \frac{y}{x} \right) \log N_{S}
    \right\}
    \notag \\
    & - c_{15, 3} \left\{
    n_{K} \frac{\left( \log \frac{y}{x} \right) (\log y)}{x \log x}
    \right\}
    \notag \\
    & - c'_{6} \frac{x^{-2}}{n} \log d_{L}
    \notag
\end{align}
where $c_{15, v} \left\{ \ldots \right\}$ comes from $S_{1, v}$ and
$c'_{6} \left\{ \ldots \right\}$ comes from $T_{1}$.

\medskip

Fix any positive constant $\epsilon$, and set $y = x^{1 +
\epsilon}$, $x = {(d_{L} N_{S}^{n})}^{C} $ for sufficiently large $C$,
one gets that the first term dominates over the other
terms by a large constant factor. Let us check this.

\medskip

The $c_{14} \{ \ldots \}$ term is bounded by some multiple of
\[
    \frac{1}{n} {(\log (d_{L} N_{S}^{n}))}^{2} 4^{- C c_{12}}
    = o \left( \frac{{(\log x)}^{2}}{n} \right)
\]
as $C$ goes to $\infty$, thus it is dominated over by
$\frac{1}{n} {\left( \log \frac{y}{x} \right)}^{2}$ by a large
constant factor. Also, this term is bounded by some multiple of
\[
    \frac{1}{n} {(\log d_{L})}^{2} 4^{1 - c_{12} C}
    \cdot (1 - \beta_{0}) \log d_{L}
\]
which is $o(\frac{(\log x)^{3}}{n} (1 - \beta_{0}))$ as $C$ goes to $\infty$,
thus it is dominated
over by \linebreak
$\frac{1}{n} {\left( \log \frac{y}{x} \right)}^{3} (1 - \beta_{0})$
by a large constant factor. From the discussion above one can verify
this assertion for the $c_{14} \{ \ldots \}$ term.

\medskip

Since
\[
    \frac{1}{x^{2}} \log \frac{y}{x}\log N_{S} <<
     \frac{1}{d_{L}^{2 C}}
    \log \frac{y}{x}
    << \frac{1}{n} \log \frac{y}{x}
\]
thus one can verify this assertion for the $c_{15, 2} \{ \ldots \}$
term.

\medskip

Other terms are easy to check. So one draws a contradiction, and we
get Theorem ~\ref{TM:A} in this case.

\medskip

Furthermore, we consider the case
\[
    1 - \beta_{0} \leq {c_{7}}^{2} {(\log d_{L} 3^{n_{L}})}^{-2},
\]
where we will use the second kernel function $k_{2} (s)$. In this
case
\[
    \log \frac{c_{7}}{(1 - \beta_{0}) \log (d_{L} 3^{n_{L}} N_{S}^{n})}
    \ge \frac{1}{2} \log {(1 - \beta_{0})}^{-1}
\]

If $\rho = \beta + i \gamma$ is a zero of $\zeta_{L} (s)$ with
$|\gamma| \leq 1$, and $\rho \ne \beta_{0}$, then by the
Deuring--Heilbronn,
\begin{align}
    |k_{2} (\rho)| &<< x^{\beta^{2} + \beta} << x^{1 + \beta}
    \notag \\
    &= x^{2} \exp \left\{ - c_{19} \frac{\log x \log {(1 - \beta_{0})}^{-1}}
    {\log (d_{L} N_{S}^{n})}\right\}
    x^{2} {(1 - \beta_{0})}^{c_{19} \frac{\log x}{\log (d_{L} N_{S}^{n})}}
    \notag
\end{align}
for some positive absolute constant $c_{19}$. Thus
\[
    \sum_{|\gamma| \leq 1 , \rho \ne \beta_{0}} |k_{2} (\rho)| <<
    x^{2} {(1 - \beta_{0})}^{c_{19} \frac{\log x}{\log (d_{L} N_{S}^{n})}} \log (d_{L} N_{S}^{n})
\]

\medskip

If $\rho = \beta + i \gamma$ is a zero of $\zeta_{L} (s)$ with
$|\gamma| \ge 1$, and $\rho \ne \beta_{0}$, we have
\[
    |k_{2} (\rho)| \leq x^{2 - \gamma^{2}} << x
\]
Thus assume $x > 2$.  Applying Lemma ~\ref{T:301}, we have
\begin{align}
    \sum_{|\gamma| > 1} |k_{2} (\rho)|
    &<< \sum_{n \ge 1} N (|2 n|)\, x^{1 + 4 n - 4 n^{2}}
    \notag \\
    &<< x \log d_{L} \sum_{n \ge 1} 2^{4 n - 4 n^{2}}
        + x \, n_{L} \sum_{n \ge 1}
    2^{4 n - 4 n^{2}} \log (2 n + 1)
    \notag \\
    &<< x \log d_{L}
    \notag
\end{align}
where $N (T)$ is the number of zeros of $\zeta_{L} (s)$ in the
region $[0, 1] \times [T - 1, T + 1]$.

\medskip

Thus
\begin{align}
    \label{4A-2} \tag{4A-2}
    k_{2} (1) & - \sum_{\rho} k_{2} (\rho)
    \ge \frac{x^{2}}{10} \min \left\{ 1, (1 - \beta_{0}) \log x \right\}
    \\
    & - c_{20} x \log d_{L} - c_{21}
    x^{2} {(1 - \beta_{0})}^{c_{19} \frac{\log x}{\log (d_{L} N_{S}^{n})}}
    \cdot \log (d_{L} N_{S}^{n})
    \notag
\end{align}
for some absolute positive constants $c_{20}$ and $c_{21}$.

\medskip

We now complete the proof of Theorem ~\ref{TM:A} in the case
\[
    1 - \beta_{0} \leq {c_{7}}^{2} {(\log (d_{L} 3^{n_{L}} N_{S}^{n}))}^{2}.
\]

Assume that for any $\fP$ in $P_{1} (C, S)$, $N \fP > x^{3 +
\delta}$. Then
\begin{align}
    0 &= \tilde{I}_{2} = \sum_{\fP \in P_{1}(C, S), N \fP \leq x^{3 + \delta}}
    (\log N\fP) \hat{k}_{2} (N \fP)
    \notag \\
    & \ge \frac{x^{2}}{10 n}
    \min \left\{ \,1, (1 - \beta_{0}) \log x \,\right\}
    \notag \\
    & - c_{20} \frac{1}{n} \log d_{L}
    \notag \\
    & - c_{21} \frac{1}{n} x^{2} {(1 - \beta_{0})}^{c_{19} \frac{\log x}{\log (d_{L} N_{S}^{n})}}
    \cdot \log (d_{L} N_{S}^{n})
    \notag
\end{align}
\begin{align}
    & - c_{22, 1} \left\{
    \frac{1}{n} {(\log x)}^{\frac{1}{2}} \log d_{L}
    \right\}
    \notag \\
    & - c_{22, 2} \left\{
    {(\log x)}^{\frac{1}{2}} \log N_{S}
    \right\}
    \notag \\
    & - c_{22, 3} \left\{
    n_{K} x^{\frac{7}{4}}
    \right\}
    \notag \\
    & - c_{22, 4} \left\{
    n_{K} x^{2 - \frac{\delta^{2}}{4}} (\log x)
    \right\}
    \notag \\
    & - c'_{6} \log d_{L}
    \notag
\end{align}
where $c_{22, v} \left\{ \ldots \right\}$ comes from $S_{2, v}$ and
$c'_{6} \left\{ \ldots \right\}$ comes from $T_{1}$.

\medskip

Fix any positive constant $\epsilon'$, and set $x = {d_{L}}^{C}$ for
sufficiently large $C$. One gets that the first term
dominates over the other terms by a large constant factor. Let us
check this.

\medskip

First be aware that by the Deuring--Heilbronn (Corollary
~\ref{T:305}), and the fact that $d_{L}^{\varepsilon} >> \log d_{L}
>> n_{L}$ for any $\varepsilon > 0$, the first term dominates over
$n_{K} x^{2 - \alpha}$ for $C$ sufficiently large for
any $\alpha > 0$.

\medskip

The $c_{21} \{ \ldots \}$ term is bounded by some multiple of
\[
    \frac{1}{n} x^{2} {(1 - \beta_{0})}^{c_{19} C} \log (d_{L} N_{S}^{n})
\]
as $C$ goes to $\infty$, and thus it is dominated over by
the first term.

\medskip

Since
\[
    {(\log x)}^{2} \log N_{S}<<
    x^{\frac{2}{1 + \epsilon'}}
\]
thus one can verify this assertion for the $c_{22, 2} \{ \ldots \}$
term.

\medskip

Other terms are easy to check now. So one draws a contradiction, and
we prove Theorem ~\ref{TM:A} in this case.

\medskip

The proof of Theorem ~\ref{TM:B} is also similar with slight modification.
Just consider $\zeta_{K} (s) L (s, \chi)$ instead of $\zeta_{K} (s)$
and we can imitate \cite{L-M-O79} to prove, and moreover use $d_{K} N (\chi) N_{S}$
instead of $d_{K} N (\chi)$

\bigskip

\section{Landau Method} \label{S:4}

In this section, we'll prove Theorem ~\ref{T:401} (Theorem ~\ref{TM:C}). We are using Landau's
idea (\cite{Landau27}), and the proof also follows \cite{LW2009}.

\medskip

\begin{theorem}  \label{T:401}
Let $\pi$ and $\pi'$ be two unitary cuspidal automorphic representations
of $\GL_{d} (K)$. Let $S$ be a finite set of places of $K$,
and $Q = {\rm max} (C (\pi), C (\pi'))$ and assume that the bound for Ramanujuan for
$\pi$ and $\pi'$ are $< R$.

\medskip

Then if $\pi \not\cong \pi'$,
there exists a place $v$ of $K$ such that $\pi_{v} \not\cong \pi'_{v}$
and
\[
N (\fP_{v}) \leq \begin{cases}
C Q^{1 + \epsilon} N_{S}^{\epsilon} &(d = 1) \notag \\
C Q^{2 d + \frac{d (d - 2)}{d H + 1} + \epsilon}
    N_{S}^{\frac{d^{3} (2 R + H)}{d H + 1} + \epsilon}
&(\text{general $d$}) \notag
\end{cases}
\]
where $C$ is some effectively computable constant
only depending on arbitrarily chosen number $H > 2 R, \epsilon > 0$, $K$ and $d$.
\end{theorem}

\medskip

\emph{Proof of Theorem ~\ref{T:401} and ~\ref{TM:C}}:

\medskip

Before we start, we pose a condition on $\pi$ and $\pi'$: Let $\delta < 1/2$
be any given positive number, and we write the infinite part of
$L (s, \tilde{\pi} \times \pi')$ as:
$\prod_{j = 1}^{d^{2} n_{K}} \Gamma_{\bR} (s + b_{j} (\tilde{\pi} \times \pi'))$.

\medskip

\textbf{(AA-Additional Assumption)}:
The horizontal distance from $\pm H$ and $0$ and all $b_{j} (\tilde{\pi} \times \pi')$
inside $\bC / \bZ$ are all greater than $\delta$, namely, for each integer $N$,
$|H - N| < \delta$ and $|\pm H - N - {\rm Re} b_{j} (\tilde{\pi} \times \pi')| < \delta$.

\medskip

We want to prove theorem ~\ref{T:401} with the assumption (AA) first.

\medskip

Form
\[
S (X, \tilde{\pi} \times \pi', S) = \sum_{n = 1}^{\infty} a_{\tilde{\pi} \times \pi', S} (n)
 \omega (\frac{n}{X})
\]
where $a_{\tilde{\pi} \times \pi', S} (n)$ is the $n$-th coefficient of
the incomplete $L$-function
\[
L^{S} (s, \tilde{\pi} \times \pi') =
\prod_{v \notin S} L (s, \tilde{\pi_{v}} \times \pi'_{v}) =
\sum_{n = 1}^{\infty} a_{\tilde{\pi} \times \pi', S} (n) n^{-s}
\]
and the weight function $\omega (x)$ defined as a smooth function which
may be specified as \cite{Wyh2006}:
\[
\omega (X) = \begin{cases}
0 &\text{($x \leq 0$ or $x \geq 3$)} \\
e^{-1/x} &\text{($0 < x \leq 1$)} \\
e^{-1/(3 - x)} &\text{($2 \leq x < 3$)} \\
\leq 1 &\text{(all $x$)}
\end{cases}
\]

\medskip

Consider the Mellin transform
\[
W (s) = \int_{0}^{\infty} \omega (x) x^{s - 1} d\,x
\]
which is a analytic function of $s$. Fix $\sigma < 0$ and let $s = \sigma + i t$
then
\[
W (s) <<_{A, \sigma} \frac{1}{{(1 + |t|)}^{A}}
\]
for all $A > 0$ by repeated partial integration.
By Mellin inversion,
\[
\omega (x) = \frac{1}{2 \pi i} \int_{(2)} W (s) x^{-s} d\, s
\]
where the integration is made along the vertical line ${\rm Re} s = 2$.

\medskip

Then we have
\begin{align}
S (X, \tilde{\pi} \times \pi', S)
 &= \frac{1}{2 \pi i} \sum_{n = 1}^{\infty} a_{\tilde{\pi} \times \pi', S} (n) \int_{(2)}
  W (s) {\left( \frac{n}{X} \right)}^{-s} d\, s \notag \\
&= \frac{1}{2 \pi i} \int_{(2)} X^{s} W (s) L^{S} (s, \tilde{\pi} \times \pi') \notag
\end{align}
where the interchange of the summation and the integral is guaranteed by the absolute
convergence along the real line $\sigma = 2$. By the standard arguments,
plus the fact that all incomplete $L$-functions of nontrivial characters are entire
of order $1$,
we may shift the integral line to get
\[
S (X, \tilde{\pi} \times \pi', S) = \frac{1}{2 \pi i} \int_{-H}
 x^{s} W (s) L^{S} (s, \tilde{\pi} \times \pi') d\,s
\]
where $H > 0$ is to be specified later.

\medskip

Let $L_{S} (s, \tilde{\pi} \times \pi') = \prod_{v \in S} L (s, \tilde{\pi}_{v} \times \pi'_{v})$,
and $L_{\infty} (s, \tilde{\pi} \times \pi') = L (s, \tilde{\pi}_{\infty} \times \pi'_{\infty})$ the
gamma factor of the Rankin-Selberg product $L$-function. Then we have the following
functional equation
\[
L^{S} (s, \tilde{\pi} \times \pi') = W (\tilde{\pi} \times \pi')
 {N (\tilde{\pi} \times \pi')}^{(1/2 - s)} G_{0} (s) G_{1} (s)
  L^{S} (1 - s, \pi \times \tilde{\pi}')
\]
where $W (\tilde{\pi} \times \pi')$
 is the root number of $\chi$ which has absolute value $1$,
\[
G_{0} (s) = \frac{L_{\infty} (1 - s, \pi \times \tilde{\pi}')}
{L_{\infty} (s, \tilde{\pi} \times \pi')}
\]
and
\[
G_{1} (s) = \frac{L_{S} (1 - s, \pi \times \tilde{\pi}')}
{L_{S} (s, \tilde{\pi} \times \pi')}
\]
We need to estimate $G_{0} (s)$ and $G_{1} (s)$ along the vertical line
$\sigma = -H$, avoiding to the pole of them.

\medskip

To be continued.

\medskip

\begin{lemma} \label{T:402}

\medskip

Under the assumption (AA),

\medskip

\textnormal{(1)}
then
\[
G_{0} (-H + iT) <<_{H, d, K, \delta} {( 1 + |t|)}^{n_{K} d^{2} (1/2 + H)}
\prod_{j = 1}^{d^{2} n_{K}} (1 + |b_{j} (\tilde{\pi} \times \pi')|)
\]

\medskip

\textnormal{(2)}
$G_{1} (-H + iT) <<_{H, d, K, \delta, \epsilon'} N_{S}^{d^{2}(2 R + H) + \epsilon'}$.
\end{lemma}

\medskip

We quote the following results on gamma function for which the proof
can be found in a lot of analysis textbooks.

\medskip

\begin{lemma}  \label{T:403}
\textnormal{\textbf{(Stirling formula)}}

\[
|\Gamma (\sigma + it)| = \sqrt{2 \pi} e^{- \frac{\pi}{2} t} {|t|}^{\sigma - 1/2}
    \left( 1 + O_{\sigma, \delta, A} \left( \frac{1}{1 + |t|} \right) \right)
\]
for all $|t| > A > 0 $ and $s = \sigma + it$ away from any poles by at least distance $\delta$.
\end{lemma}

\medskip

\qedsymbol

\medskip

\emph{Proof of Lemma ~\ref{T:402}}:

\medskip

(1) Write $b_{j} = b_{j} (\tilde{\pi} \times \pi') = u_{j} + i v_{j}$,
and we have $b_{j} (\pi \times \tilde{\pi'})
= \overline{b_{j} (\tilde{\pi} \times \pi')} = \overline{b_{j}} = u_{j} - i v_{j}$.
Put $s = \sigma + i t$.

\medskip

\begin{align}
G_{0} (s) &= \frac{L_{\infty} (1 - s, \pi \times \tilde{\pi}')}
{L_{\infty} (s, \tilde{\pi} \times \pi')} \notag\\
&= \pi^{-d^{2} n_{K} / 2 + d^{2} n_{K} s} \prod_{j = 1}^{d^{2} n_{K}}
\frac{\Gamma ((1 - s + \overline{b_{j}}) / 2)}{\Gamma ((s - b_{j}) / 2)} \notag \\
&<<_{\sigma, \delta, d, K} \prod_{j = 1}^{d^{2} n_{K}}
 \frac{{|t + v_{j}|}^{\frac{1 - \sigma + v_{j}}{2} - \frac{1}{2}}}
 {{|t + v_{j}|}^{\frac{\sigma + v_{j}}{2} - \frac{1}{2}}} \notag \\
&<<_{\sigma, \delta, d, K} \prod_{j = 1}^{d^{2} n_{K}} {|t + v_{j}|}^{1/2 - \sigma} \notag
\end{align}

\medskip

Hence, under the assumption (AA)
\[
G_{0} (-H + i t) <<_{\sigma, d, K, H, \delta}
{(1 + |t|)}^{d^{2} n_{K} (1/2 + H)} \prod_{j = 1}^{d^{2} n_{K}}
{(1 + |b_{j}|)}^{1/2 + H}
\]

(2) For each $v \in S$, write
$L (s, \pi_{v}) = \prod_{j = 1}^{d} (1 - a_{v, j} q_{v}^{-s})$
and
$L (s, \pi'_{v}) = \prod_{j = 1}^{d} (1 - a'_{v, j} q_{v}^{-s})$.
Then
$L (s, \tilde{\pi}_{v} \times \pi'_{v}) =
\prod_{j, k = 1 \ldots d} (1 - \overline{a_{v, j}} a'_{v, k} q_{v}^{-s})$
and hence
$L (s, \pi_{v} \times \tilde{\pi}'_{v}) =
\prod_{j, k = 1 \ldots d} (1 - a_{v, j} \overline{a'_{v, k}} q_{v}^{-s})$

\medskip

\begin{align}
G_{1} (s) &= \frac{L_{S} (1 - s, \pi \times \tilde{\pi}')}
{L_{S} (s, \tilde{\pi} \times \pi')} \notag\\
&= \prod_{v \in S} \frac{L (1 - s, \pi_{v} \times \tilde{\pi}'_{v})}
{L (s, \tilde{\pi}_{v} \times \pi'_{v})} \notag\\
&= \prod_{v \in S} \prod_{j, k = 1 \ldots d}
\frac{1 - \overline{a_{v, j}} a'_{v, k} q_{v}^{-s}}{1 - a_{v, j} \overline{a'_{v, j}} q_{v}^{s - 1}} \notag
\end{align}

Note that by and by the assumption on the Ramanujuan bounds and the results on it
(\cite{LRS99}, \cite{BB}), we have:
$|a_{v, j}|, |a'_{v, k}| < q_{v}^{R} < \leq q_{v}^{1/2 - 1/ (d^{2} + 1)}$.
Hence,
\begin{align}
G_{1} (- H + i t) &\leq \prod_{v \in S} {\left( \frac{1 + q_{v}^{2 R + H}}
    {1 - q_{v}^ {- (H + 2 / (d^{2} + 1))}} \right)}^{d^{2}} \notag \\
&\leq \prod_{v \in S} {\left( \frac{2 q_{v}^{2 R + H}}
    {1 - 2^{- (H + 2 / (d^{2} + 1))}} \right)}^{d^{2}} \notag \\
&\leq 2^{|S| d^{2}} N_{S}^{d^{2}(2 R + H)}
{1 - 2^{- (H + 2 / (d^{2} + 1))}}^{-|S|}
\notag \\
&<<_{d, K, \epsilon'} N_{S}^{d^{2} (2 R + H) + \epsilon'}
\notag
\end{align}
Here we use the estimation $a^{|S|} <<_{K, d, \epsilon'} N_{S}^{\epsilon'}$.

\medskip

\qedsymbol

\medskip

\emph{Proof of Theorem ~\ref{T:401}, Cont.\ }:

\medskip

Still under the asumption (AA).

Now
\begin{align}
&S (X, \tilde{\pi} \times \pi', S) \notag \\
&=
\frac{1}{2 \pi i} \int_{(-H)} X^{s} W (s) W (\tilde{\pi} \times \pi')
 {A (\tilde{\pi} \times \pi')}^{1/2 - s} G_{0} (s) G_{1} (s)
        L^{S} (1 - s, \pi \times \tilde{\pi}') d\,s \notag \\
&=
\frac{1}{2 \pi i} \int_{(-H)} X^{s} W (s) W (\tilde{\pi} \times \pi')
 {A (\tilde{\pi} \times \pi')}^{1/2 - s} G_{0} (s) G_{1} (s)
        \sum_{n = 1}^{\infty} \frac{a_{\tilde{\pi} \times \pi', S} (n)}{n^{1 + H}}  d\,s \notag \\
&=
\frac{1}{2 \pi i} \sum_{n = 1}^{\infty}
 \frac{a_{\tilde{\pi} \times \pi', S} (n)}{n^{1 + H}}
        \int_{(-H)} X^{s} W (s) W (\tilde{\pi} \times \pi')
         {N (\tilde{\pi} \times \pi')}^{1/2 - s} G_{0} (s) G_{1} (s)
        n^{s + H} d\,s \notag
\end{align}

Here the interchange of the sum and the integral is guaranteed by the absolute convergence
of the Dirichlet series, rapid decay of $W (s)$.

\medskip

Note that when $H > 2 R$, we have
\begin{align}
\sum_{n = 1}^{\infty} & |\frac{a_{\tilde{\pi} \times \pi', S} (n)}{n^{1 + H}}| \notag \\
&\leq \sum_{n = 1}^{\infty} \zeta_{K} (1 + H - 2 R) \notag
\end{align}
by the arguments in the preliminary subsection on Bound for Ramanujuan.

\medskip

Applying Lemma ~\ref{T:402}, and the estimation above,
we have
\begin{align}
&S (X, \tilde{\pi} \times \pi', S) \notag \\
&<<_{H, K} \zeta_{K} (1 + H - 2 R)
\int_{(-H)} \left| X^{s} W (s) W (\tilde{\pi} \times \pi')
     {N (\tilde{\pi} \times \pi')}^{1/2 - s} G_{0} (s) G_{1} (s)
         d\,s \right| \notag \\
&<<_{H, K, \delta, \epsilon}
X^{-H} {N (\tilde{\pi} \times \pi')}^{1/2 + H}
\notag \\ &\qquad
\prod_{j = 1}^{d^{2} n_{K}} {(1 + |b_{j} (\tilde{\pi} \times \pi')|)}^{1/2 + H}
N_{S}^{d^{2} (2 R + H) + \epsilon} \int_{(-H)}
    W (s) {(1 + |t|)}^{(1/2 + H) n_{K}} d\,s \notag \\
&<<_{H, K, \delta, \epsilon}
 X^{-H} {C (\tilde{\pi} \times \pi')}^{1/2 + H} N_{S}^{d^{2} (2 R + H) + \epsilon}
 \notag
\end{align}

\medskip

To establish the theorem, we need to bound $S (X, \tilde{\pi}, \times \pi', S)$
below. Now assume that $\pi_{v} \cong \pi'_{v}$ for all $v \notin S$ such that $N (\fP_{v}) \leq 3 X$.
Then
\[
S (X, \tilde{\pi} \times \pi', S) = S (X, \tilde{\pi} \times \pi, S)
\]

\medskip

Now,
we have
\[
S (X, \tilde{\pi} \times \pi, S) >>_{K, d} X^{1/d} / (\log (X)) - |S|
\]
since by the prime number theory, when $X$ is large,
there is a prime $\fP_{v}$ of $K$ such that $X < q_{v}^{d^{2}} < 2 X$
(prime number theorem and Bertrand)
where $q_{v} = N \fP_{v}$.

\medskip

Thus, by Proposition ~\ref{T:20X},
$S (X, \tilde{\pi} \times \pi, S)$ is greater than $e^{-1}$
multiples of the primes $\fP_{v}$ of $K$
of degree $1$ outside $S$ such that $X \leq \fP_{v} \leq 2 X$,
when $X > 4 A {|S|}^{2}$, $S (X, 1, S) >>_{K, d} X^{1/d} / (\log (X))$
for large $A > 0$.

\medskip

Then we have, when $X > 4 A {|S|}^{2}$,
\begin{align}
X^{1/d} / (\log (X)) &<<_{K, d, \delta} S (X, \tilde{\pi} \times \pi, S)
 = S (X, \tilde{\pi} \times \pi', S) \notag \\
 < C' X^{-H} {C (\tilde{\pi} \times \pi')}^{1/2 + H}
  N_{S}^{d^{2} (2 R + H) + \epsilon}
\notag
\end{align}
for some $C'$ depending on $H$, $K$, $d$, $\delta$.

\medskip

Therefore,
\begin{align}
X^{H + 1/d - \epsilon} &<<_{\epsilon} X^{H + 1/d} {\log (X)}^{-1} \notag \\
&<<_{H, K, d, \epsilon, \delta} {C (\tilde{\pi} \times \pi')}^{1/2 + H} N_{S}^{d^{2} (2 R + H) + \epsilon} \notag
\end{align}
and thus
\[
X <<_{\epsilon, H, K, d, \delta} {C (\tilde{\pi} \times \pi')}^{\frac{H + 1/2}{H + 1/d - \epsilon}}
 N_{S}^{\frac{d^{2} (2 R + H)}{H + 1/d - \epsilon}}
\]

\medskip

\[
X < C Q^{2 d + \frac{d^{2} (d - 2)}{d H + 1} + \epsilon} N_{S}^{\frac{d^{3} (2 R + H)}{d H + 1} + \epsilon}
\]
where $C$ depends on $\epsilon$ and $K$, $d$, $\delta$ and $H > 2 R$.
As
\[
X < 4 A {|S|}^{2}
<<_{\epsilon, A, K} N_{S}^{\epsilon}
\]

\medskip

Hence Theorem ~\ref{T:401} follows under the assumption (AA).

\medskip

Now in general, we can choose small $\delta$ depending only on $K$,
$d$ and $\epsilon$ such that the assumption (AA) holds always for some $H'$
with $|H' - H| < \delta /2$.  Then we still get the theorem
for $H'$ which is sufficiently close to $H$.

\medskip

\qedsymbol

\end{document}